\documentclass[review,11pt]{elsarticle}
\usepackage{algorithm}
\usepackage{algpseudocode}
\usepackage{amsmath}
\usepackage{amssymb}
\usepackage{bm}
\usepackage{booktabs}
\usepackage{color}
\usepackage{epsfig}
\usepackage{epstopdf}
\usepackage{flafter}
\usepackage{float}
\usepackage{graphicx}
\usepackage{geometry}
\usepackage[colorlinks,linkcolor=red,anchorcolor=gray,citecolor=green,urlcolor=black]{hyperref}
\usepackage{ifthen}
\usepackage{lineno}
\usepackage{multirow}
\usepackage{makecell} 
\usepackage{pifont}

\usepackage{subfigure}
\usepackage{times}
\usepackage{threeparttable}
\usepackage{url}  
\biboptions{numbers,sort&compress}
\modulolinenumbers[5]
\journal{Journal of Computational and Applied Mathematics}

\newtheorem{thm}{\textbf{Theorem}}[section]
\newtheorem{mydef}{\textbf{Definition}}[section]
\newtheorem{lem}{\textbf{Lemma}}[section]

\begin{document}
\begin{frontmatter}
\title{Multi-dimensional imaging data recovery via minimizing the partial sum of tubal nuclear norm}

\author[swufe]{Tai-Xiang Jiang}
\ead{taixiangjiang@gmail.com}

\author[uestc]{Ting-Zhu Huang\corref{cor}}
\ead{tingzhuhuang@126.com.}

\author[uestc]{Xi-Le Zhao\corref{cor}}
\ead{xlzhao122003@163.com.}

\author[uestc]{Liang-Jian Deng}
\ead{liangjiand1987112@126.com}

\cortext[cor]{Corresponding authors.}

\address[uestc]{School of Mathematical Sciences/Research Center for Image and Vision Computing, University of Electronic Science and Technology of China, Chengdu, Sichuan, 611731, China}
\address[swufe]{FinTech Innovation Center, Financial Intelligence and Financial Engineering Research Key Laboratory of Sichuan province, School of Economic Information Engineering, Southwestern University of Finance and Economics, Chengdu, Sichuan, 611130, China}
\begin{abstract}
In this paper, we investigate tensor recovery problems within the tensor singular value decomposition (t-SVD) framework. We propose the partial sum of the tubal nuclear norm (PSTNN) of a tensor. The PSTNN is a surrogate of the tensor tubal multi-rank. We build two PSTNN-based minimization models for two typical tensor recovery problems, i.e., the tensor completion and the tensor principal component analysis. We give two algorithms based on the alternating direction method of multipliers (ADMM) to solve proposed PSTNN-based tensor recovery models. Experimental results on the synthetic data and real-world data reveal the superior of the proposed PSTNN.
\end{abstract}
\end{frontmatter}
\section{Introduction}\label{sec:Intro}
The tensor, a multi-dimensional extension of the matrix, is an important data format and has been applied in lots of real-world applications, for example, the video data recovery \cite{jiang2018fastderain}, the hyperspectral data recovery and fusion \cite{zhuang2018fast,li2018fusing}, the personalized web search \cite{Sun2005web}, and seismic data reconstruction \cite{Kreimer2012HSVDtensor}. Among these applications, how to accurately characterize and rationally utilize the inner structure of these multi-dimensional data is of crucial importance \cite{Liu2013PAMItensor}.

In the matrix processing, low-rank models can robustly and efficiently handle two-dimensional data of various sources, and the solutions are generally theoretically guaranteed in many applications \cite{candes2012exact,ma2017truncated}.
However, how to extend the low-rank definition from matrices to tensors is still an open problem.
The most two popular tensor rank definitions in the past decade are the CANDECOMP/PARAFAC (CP)-rank, which is related to the CANDECOMP/PARAFAC decomposition \cite{kiers2000towards}, and Tucker-rank (or denoted as ``$n$-rank'' in \cite{gandy2011tensor}), which is corresponding to the Tucker decomposition \cite{tucker1966some,kapteyn1986approach}.

In this paper, we fix our attention on a newly emerged tensor decomposition paradigm, the tensor singular value decomposition (t-SVD), and the notion of the tensor rank derived from t-SVD, i.e., the tubal multi-rank.
The t-SVD was been initially proposed in \cite{braman2010third, kilmer2011factorization} and it allows new extensions of familiar matrix analysis to the tensor while avoiding the loss of information inherent in matricization or flattening of the tensor \cite{kilmer2013third}.
The tubal nuclear norm, which is a convex surrogate of the tubal multi-rank, is utilized to handle the tensor completion problem by Zhang {\it et al.} \cite{zhang2017exact} and the tensor completion from sparsely corrupted observations by Jiang {\it et al.} \cite{jiang2017exact}.

The t-SVD is defined based on the tensor-tensor product (t-prod). Owing to its particular structure, the t-prod is equivalent to the matrix-matrix product after the Fourier transform. Meanwhile, according to the definition, the TNN is equivalent to the matrix nuclear norm of the block diagonal unfolding of the Fourier transformed tensor (See Eq. \eqref{tnn} in the Appendix for details).
However, in the matrix case, minimizing the nuclear norm would cause some unavoidable biases \cite{oh2013partial, oh2016partial}.
For example, the variance of the estimated data would be smaller than the original data when equally shrink every singular value. Similarly, the estimated results may be lower-rank than the original data. Therefore, following the research path in \cite{gu2014weighted,lu2015generalized,oh2016partial}, we consider minimizing the proposed partial sum of the tubal nuclear norm (PSTNN), which only consists of the small singular values.
On the one hand, minimizing the PSTNN would directly shrink the small singular values without any actions on the large ones, resulting in low tubal multi-rank estimations without rank deficiency situations. On the other hand, the corresponding minimization problem is easy to optimize with the proposed solver.

The main contributions are the following three aspects.
First, we propose a surrogate of the tensor tubal multi-rank, i.e., the PSTNN.
Second, to optimize the PSTNN-based minimization Problem, we extend the partial singular value thresholding (PSVT) operator, which was primarily proposed in \cite{oh2013partial}, for the matrices in the complex field, and demonstrate that it is the exact solution to the PSTNN-based minimization problem.
Third, we propose two PSTNN-based models to solve the typical tensor recovery problems, i.e., the tensor completion and the tensor robust principal component analysis. Afterward, two alternating direction method of multipliers (ADMM) algorithms using the PSVT solver are developed to optimize two PSTNN based models.
Moreover, we conduct experiments on synthetic data and real-world data. The results illustrate that proposed methods can effectively handle tensor recovery problems.

The organization of this paper is organized as follows.
Section \ref{sec:Pre} presents some preliminaries.
In Section \ref{sec:Main}, we give the main results.
Section \ref{sec:Exp} reports the experimental results.
Finally, in Section \ref{sec:Con}, some conclusions are drawn.
\section{Notation and preliminaries}\label{sec:Pre}

Before giving the main results, we briefly introduce the basic tensor notations and exhibit the t-SVD algebraic framework. The notations and definitions in this section are referred to \cite{kolda2009tensor,kilmer2013third,kilmer2011factorization,lu2018tensor,zhang2017exact}.

Throughout this paper, lowercase letters, e.g., $x$, boldface lowercase letters, e.g., $\mathbf x$, boldface upper-case letters, e.g., $\mathbf X$, and boldface calligraphic letters,  e.g., $\mathbf{\mathcal{X}}$, are respectively used to denote scalars, vectors, matrices, and tensors.
The $(i_{1},i_2,\cdots,i_{N})$-th  element of an $N$-mode tensor is denoted as $x_{i_{1}i_2\cdots i_{N}}$.
The inner product of two tensors $\mathcal{X}$ and $\mathcal{Y}$, of the same size, is defined as $\langle\mathcal{X},\mathcal{Y}\rangle:=\sum\limits_{i_{1},i_{2},\cdots,i_{N}}x_{i_{1}i_2\cdots i_{N}}\cdot y_{i_{1}i_2\cdots i_{N}}$.
Then, the tensor Frobenius norm of $\mathcal{X}$ is defined as $\left\|\mathcal{X}\right\|_{F}:=\sqrt{\langle\mathcal{X},\mathcal{X}\rangle} = \sqrt{\sum_{i_1,\cdots,i_N}x_{i_{1}i_2\cdots i_{N}}^2}$.

We denote the Fourier transform along the third mode of a third-order tensor $\mathbf{\mathcal{X}}\in \mathbb{R}^{n_{1}\times n_2\times n_{3}}$ as $\widehat{\mathbf{\mathcal{X}}}={\tt{fft}}(\mathbf{\mathcal{X}},[],3)$.
Meanwhile, the inverse transformation is denoted as $\mathbf{\mathcal{X}}={\tt{ifft}}(\widehat{\mathbf{\mathcal{X}}},[],3)$.
As shown in \cite{lu2018tensor}, $\|{\mathcal{X}}\|_{F}  = \frac{1}{\sqrt{n_3}}\|\widehat{\mathcal{X}}\|_{F}$.
To save space, the definitions related to the t-SVD framework are given in Appendix \ref{tsvddefs}. We list all the notations in Table \ref{notationstable}.

\begin{table*}[htbp]\label{notationstable}
\renewcommand\arraystretch{1.0}
\caption{Tensor notations}
 \begin{tabular}{c p{0.79\columnwidth}}
  \toprule
Notation &  Explanation \\
  \midrule
$\mathbf{\mathcal{X}},\mathbf{X},\mathbf{x},x$
                                & Tensor, matrix, vector, scalar.\\
$x_{i_{1}i_2\cdots i_{N}}$
                                & The $(i_{1},i_2,\cdots,i_{N})$-th  element of an $N$-mode tensor $\mathcal{X}$.\\
$\langle\mathbf{\mathcal{X}},\mathbf{\mathcal{Y}}\rangle$
                                & The \textbf{inner product} of two same-sized tensors $\mathbf{\mathcal{X}}$ and $\mathbf{\mathcal{Y}}$.\\
$\left\|\mathbf{\mathcal{X}}\right\|_{F}$
                                & The \textbf{Frobenius norm} of a tensor $\mathbf{\mathcal{X}}$.\\
$\widehat{\mathbf{\mathcal{A}}}$
                                & The Fourier transformed tensor of $\mathcal{A}$. \\
$\overline{\mathcal{{A}}}$
                                & The block-diagonal form unfolding of $\widehat{\mathbf{\mathcal{A}}}$. (Definition \ref{Def:bldg})\\
$\text{rank}_r(\mathbf{\mathcal{A}})$
                                & The tubal multi-rank of a tensor $\mathcal{A}\in\mathbb{R}^{n_1\times n_2\times n_3}$. (Definition \ref{Def:tubal})\\
$\|\mathbf{\mathcal{A}}\|_{\text{TNN}}$
                                & The tubal nuclear norm (TNN) of a tensor $\mathcal{A}\in\mathbb{R}^{n_1\times n_2\times n_3}$. (Definition \ref{Def:TNN})\\

\bottomrule
 \end{tabular}
 \end{table*}

\section{Main results}\label{sec:Main}
Minimizing the rank surrogate to enhance the low-dimensionality of the underlying target data is an effective way to recover the multidimensional imaging data, which is naturally in the tensor format, from incomplete or corrupted observed data. The tubal nuclear norm (TNN) is minimized to enhance the low tubal multi-rank property of the multi-dimensional visual data for the tensor completion problem in \cite{zhang2014novel,zhang2017exact}.
The tensor nuclear norm, which is similar to the TNN but defined with a factor $\frac{1}{n_3}$ in \cite{lu2018tensor}, is also minimized to promote the low-rankness for handling the RPCA problem  \cite{lu2018tensor} and the outlier-RPCA problem \cite{zhou2017outlier}.

In this section, the definition of the PSTNN is given at first. Then the PSVT-based solver for the PSTNN-based minimization problem is presented.
Subsequently, we propose the PSTNN-based tensor completion model and Tensor RPCA model and their corresponding algorithms, respectively.
\subsection{Partial sum of the tubal nuclear norm (PSTNN)}

%

\begin{figure}[!htp]
  \centering
  \includegraphics[width=0.6\textwidth]{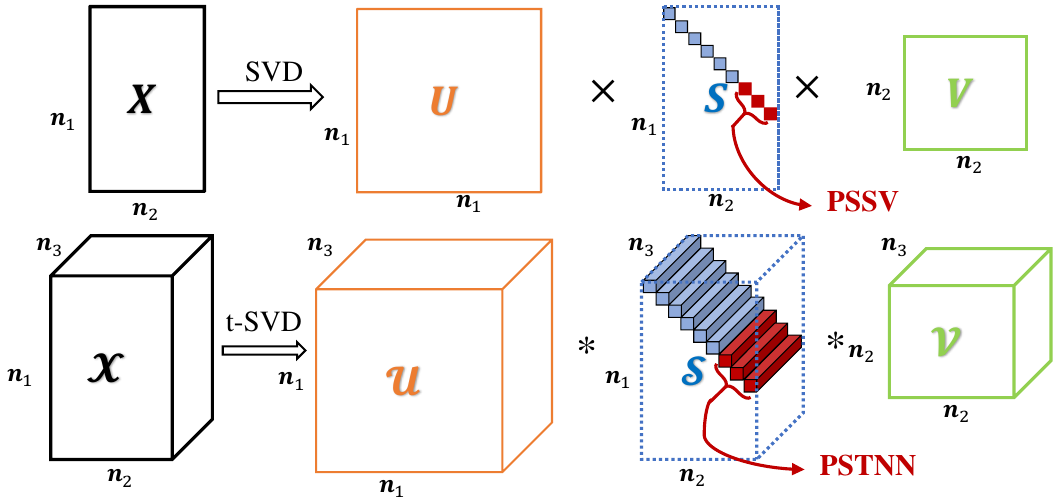}
  \caption{The illustration of the distinction and the connection between PSSV of a matrix (fisrt row) and PSTNN of a tensor (second row).}
  \label{m2t}
\end{figure}

Our PSTNN is extended from the partial sum of singular values (PSSV) \cite{oh2013partial,oh2016partial}. The PSTNN of a three way tensor $\mathbf{\mathcal{A}}\in \mathbb{R}^{n_1\times n_2\times n_3}$ is given as
\begin{equation}
\begin{aligned}
\|\mathbf{\mathcal{A}}\|_\text{PSTNN}\triangleq\sum\limits_{i=1}^{n_3}\|\widehat{\mathbf{\mathcal{A}}}^{(i)}\|_{p=N}.
\end{aligned}
\label{PTNN_def}
\end{equation}
In (\ref{PTNN_def}), $\| \cdot\|_{p=N}$ is the PSSV  \cite{oh2013partial,oh2016partial}, which is defined as $\| \mathbf{X}\|_{p=N}=\sum_{i=N+1}^{\min (m,n)} \sigma_i(\mathbf{X})$ for a matrix $\mathbf{X}\in\mathbb{C}^{n_1\times n_2}$, where $\sigma_i(\mathbf{X})$ $(i=1,\dots,\min (m,n))$ denotes its $i$-th largest singular value.
It can be observed from Figure \ref{m2t} that the proposed PSTNN is a high order extension of PSSV and the definition of PSTNN maintains an explicit meaning within the t-SVD algebraic framework, i.e., the sum of the red tubes in Figure \ref{m2t}.

It noteworthy that, 
according to Definition \ref{Def:tubal},
the $i$-th element of the tubal multi-rank of a tensor $\mathcal{A}$ is $\text{rank}( \widehat{\mathbf{\mathcal{A}}}^{ (i)} )$, and Definition \ref{Def:bldg} implies $\sum_{i=1}^{n_3}\text{rank} (\widehat{\mathbf{\mathcal{A}}}^{(i)}) = \text{rank}(\overline{\mathbf{\mathcal{A}}})$.
Thus the $l_1$ norm (sum of the absolute values of a vector) of $\mathbf{\mathcal{A}}$'s tubal multi-rank equals to the rank of its block-diagonal unfolding of $\overline{\mathbf{\mathcal{A}}}$. That is
\begin{equation}
\|\text{rank}_r(\mathcal{A})\|_1=\text{rank}(\overline{\mathbf{\mathcal{A}}}).
\end{equation}
More precisely, the TNN (defined in Definition \ref{Def:TNN}) is a convex relaxation of the $l_1$ norm of a three order tensor's tubal multi-rank. Thus, the proposed PSTNN is also a surrogate of $\|\text{rank}_r(\mathcal{A})\|_1$.

\subsection{The PSTNN-based minimization problem}

In this subsection, we introduce the general solving scheme for the PSTNN-based minimization problem, which is fundamental for solving the PSTNN-based tensor completion and PSTNN-based robust
principal component analysis problems in the subsequent two subsections.
The PSTNN-based minimization problem aims at restoring a tensor from its
observation under PSTNN regularization.
For an observed tensor $\mathbf{\mathcal{Y}}$, the PSTNN-based  minimization problem is:
\begin{equation}
\mathcal{X}^{*}=\arg\left\{\min\limits_\mathcal{X}\left(\lambda\|\mathcal{X}\|_\text{PSTNN} +\frac{1}{2}\|\mathcal{X}-\mathcal{Y}\|_F^2\right)\right\},
\label{KEY}
\end{equation}
where $\mathbf{\mathcal{X}}$ and $\mathbf{\mathcal{Y}} \in \mathbb{R}^{n_1\times n_2\times n_3}$, and $\lambda$ is non-negative parameter, which controls the balance between the PSTNN regularization and the distance to the observation.

Considering the linearity of the Fourier transform and the property that $\|\mathcal{A}\|_F^2 = \frac{1}{n_3}\|\widehat{\mathcal{A}}\|_F^2$ for any $\mathcal{X},\mathcal{Y}\in \mathbb{R}^{n_{1}\times n_2\times n_{3}}$, we have
$\frac{1}{2}\|\mathbf{\mathcal{X}}-\mathbf{\mathcal{Y}}\|_F^2 = \frac{1}{2n_3}\|\widehat{\mathcal{X}}-\widehat{\mathcal{Y}}\|_F^2
= \frac{1}{2n_3}\sum\limits_{k=1}^{n_3}\|\widehat{\mathbf{\mathcal{X}}}^{(k)}-\widehat{\mathbf{\mathcal{Y}}}^{(k)}\|_F^2$.
Meanwhile, since $\|\mathbf{\mathcal{X}}\|_\text{PSTNN}= \sum\limits_{i=1}^{n_3}\|\widehat{\mathbf{\mathcal{X}}}^{(i)}\|_{p=N}$ and $\mathcal{X} = {\tt{ifft}}(\widehat{\mathbf{\mathcal{X}}},[],3)$, the minimization problem in Eq. \eqref{KEY} is equivalent to
\begin{equation}
\begin{aligned}
\{\widehat{\mathbf{\mathcal{X}}}^{(1)*},&\widehat{\mathbf{\mathcal{X}}}^{(2)*},\cdots,\widehat{\mathbf{\mathcal{X}}}^{(n_3)*}\} \\ &=\arg\left\{\min\limits_{\widehat{\mathbf{\mathcal{X}}}^{(1)},\widehat{\mathbf{\mathcal{X}}}^{(2)},\cdots,\widehat{\mathbf{\mathcal{X}}}^{(n_3)}} \left(\lambda\sum\limits_{i=1}^{n_3}\|\widehat{\mathbf{\mathcal{X}}}^{(i)}\|_{p=N} +\frac{1}{2n_3}\sum\limits_{k=1}^{n_3}\|\widehat{\mathbf{\mathcal{X}}}^{(k)}-\widehat{\mathbf{\mathcal{Y}}}^{(k)}\|_F^2\right)\right\}.
\end{aligned}
\label{KEY2}
\end{equation}
Thus, the minimization problem in (\ref{KEY2}) can be decoupled into $n_3$ matrix minimization problems with respect to , i.e.,
\begin{equation}
\widehat{\mathbf{\mathcal{X}}}^{(k)*}=\arg\left\{\min\limits_{\widehat{\mathbf{\mathcal{X}}}^{(k)}} \left(\lambda\|\widehat{\mathbf{\mathcal{X}}}^{(k)}\|_{p=N}+\frac{\beta}{2}\|\widehat{\mathbf{\mathcal{X}}}^{(k)}-\widehat{\mathbf{\mathcal{Y}}}^{(k)}\|_F^2\right)\right\},
\label{PTNN_matrix}
\end{equation}
where  $\widehat{\mathbf{\mathcal{X}}}^{(k)},\widehat{\mathbf{\mathcal{Y}}}^{(k)}\in\mathbb{C}^{n_1\times n_2}$, $\beta = 1/n_3$, and $k=1,2,\cdots,n_3$.
The tensor optimization problem (\ref{KEY}) is herein transformed to $n_3$ matrix optimization problems in (\ref{PTNN_matrix}) in the Fourier transform domain.
It should be note that, Oh {\em et al.} have proposed the exact solution of (\ref{PTNN_matrix}), which is indeed a PSSV-based minmization problem, in \cite{oh2013partial,oh2016partial} for real matrices.
Hence, the solving scheme in  \cite{oh2013partial,oh2016partial} should be generalized to the complex matrices.


Before extending the PSVT operator for the matrices in the complex field, we first restate the von Neumann's lemma \cite{von1937some,mirsky1975trace,de1994exposed}.


\begin{lem}[von Neumann \cite{von1937some}]
If $\mathbf{A},\mathbf{B}$ are complex ${m\times n}$ matrices with singular values
\[
\sigma_1^ {A}\geq\cdots\geq\sigma_{\min (m,n)}^{A},\quad \sigma_1^{B}\geq\cdots\geq\sigma_{\min (m,n)}^{B}
\]
respectively, then
\begin{equation}
|\langle\mathbf A,\mathbf B\rangle|=|\text{Tr}(\mathbf A^{\text{H}} \mathbf B)|\leq\sum\limits_{r=1}^{\min(m,n)}\sigma_r^{A}\sigma_r^{B}.
\label{von}
\end{equation}
Moreover, equality holds in (\ref{von}) $\iff$ $\mathbf A$ and $\mathbf B$ maintains the same right and left singular vectors, i.e.,
\begin{equation}
\mathbf A=\mathbf U{\tt{diag}}(\sigma(\mathbf A))\mathbf V^{\text{H}}\  \text{and}\ \mathbf B=\mathbf U{\tt{diag}}(\sigma(\mathbf B))V^{\text{H}},
\end{equation}
\end{lem}
where $\sigma(\mathbf A)=[\sigma_1^{X},\cdots,\sigma_{\min (m,n)}^{A}]$ and $\sigma(\mathbf B)=[\sigma_1^{B}, \cdots ,\sigma_{\min (m,n)}^{B}]$.
%

Then, we restate the corresponding theorem, which utilized the von Neumann's lemma, in \cite{oh2013partial,oh2016partial} and extend it to the complex matrices case in the meantime.


\begin{thm}[PSVT]\label{PSVT_TH}
Let $\mathbf A,\mathbf B\in\mathbb{C}^{n_1\times n_2}$, which are two complex matrices, $\tau>0$, and $l=\min(n_1,n_2)$. $\mathbf B$ can be written as the linear superposition of two items, i.e., $\mathbf B=\mathbf B_1+\mathbf B_2=\mathbf U_{ B_1}\mathbf D_{B_1}\mathbf V_{B_1}^H+ \mathbf U_{ B_2}\mathbf D_{B_2}\mathbf V_{B_2}^H$, where  $\mathbf U_{B_1}$,$\mathbf V_{B_1}$ are the singular vector matrices corresponding to the $N$ largest singular values, and $\mathbf U_{B_2}$, $\mathbf V_{B_2}$ from the $(N+1)$-th to the last singular values.
A complex matrix PSSV minimization problem is
\begin{equation}\label{pssv_matrix}
\mathbf A^* = \arg\left\{\min\limits_{\mathbf A} \left(\lambda\|\mathbf A\|_{p=N} +\frac{\beta}{2}\|\mathbf A-\mathbf B\|_F^2\right)\right\}.
\end{equation}
Then, the matrix PSSV minimization in (\ref{pssv_matrix}) can be optimized by the PSVT operator as
\begin{equation}\label{psvt}
\begin{aligned}
\mathbf A^*=\mathbb{P}_{N,\tau}(\mathbf B)=\mathbf U_B(\mathbf D_{B_1}+\mathbf{\mathcal{S}}_\tau[\mathbf D_{B_2}])\mathbf V_B^H = \mathbf B_1+\mathbf U_{ B_2}\mathbf{\mathcal{S}}_\tau[\mathbf D_{B_2}])\mathbf V_{B_2}^H,
\end{aligned}
\end{equation}
where $\mathbf D_{B_1}={\tt{diag}}(\sigma_1^{B},\cdots,\sigma_N^{B},0,\cdots,0)$, $\mathbf D_{B_2}={\tt{diag}}(0,\cdots,0,\sigma_{N+1}^{B},\cdots,\sigma_l^{B},)$, and $\mathbf{\mathcal{S}}_\tau[\cdot]={\tt{sign}}(\cdot)\cdot\max(|\cdot|-\tau,0)$ $(\tau = \frac{\lambda}{\beta})$ is the soft-thresholding operator.
\end{thm}


The proof of Theorem \ref{PSVT_TH} is exhibited in Appendix \ref{ProofTH3.1}.
Then, the solution of (\ref{PTNN_matrix}) can be obtained as
\begin{equation}
{\widehat{\mathbf{\mathcal{X}}^*}}^{(k)}=\mathbb{P}_{N,\tau}\left(\widehat{\mathbf{\mathcal{Y}}}^{(k)}\right).
\end{equation}
We summarize the steps to solve (\ref{KEY}) in Algorithm \ref{KEY_alg}.




\begin{algorithm}[!htbp]
\caption{Solving (\ref{KEY}) using PSVT }
\begin{algorithmic}[1]
\renewcommand{\algorithmicrequire}{\textbf{Input:}} 
\renewcommand{\algorithmicensure}{\textbf{Output:}}
\Require
$\mathbf{\mathcal{B}}\in\mathbb{R}^{n_1\times n_2\times n_3}$, $\lambda$, the given tubal multi-rank $\text{rank}_r$
\renewcommand{\algorithmicrequire}{\textbf{Initialization:}}
\Require $\widehat{\mathbf{\mathcal{A}}}=\tt{zeros}(n_1\times n_2\times n_3)$,  $\beta=\frac{1}{n_3}$
\State $\widehat{\mathbf{\mathcal{B}}}\leftarrow\tt{fft}(\mathbf{\mathcal{B}},[],3)$, $\tau\leftarrow\frac{\lambda}{\beta}$
\For {$k=1:n_3$}
\State
$N\leftarrow$ the $k$-th element of $\text{rank}_r$
\State $\widehat{\mathbf{\mathcal{A}}}^{(k)}\leftarrow\mathbb{P}_{N,\tau}\left(\widehat{\mathbf{\mathcal{B}}}^{(k)}\right)$
\EndFor
\State $\mathbf{\mathcal{A}}\leftarrow\tt{ifft}(\widehat{\mathbf{\mathcal{A}}},[],3)$
\Ensure
$\mathbf{\mathcal{A}}\in\mathbb{R}^{n_1\times n_2\times n_3}$
\end{algorithmic}
\label{KEY_alg}
\end{algorithm}

In the following subsections, based on the proposed rank approximation, we can easily give our proposed tensor completion model and tensor RPCA model.


\subsection{Tensor completion using PSTNN}

A tensor completion model using PSTNN is given as
\begin{equation}
\begin{aligned}
\min\limits_{\mathbf{\mathcal{X}}} \quad&\|\mathbf{\mathcal{X}}\|_{\text{PSTNN}}\\
\text{s.t.}\quad& \mathcal{P}_\Omega(\mathcal{X})=\mathcal{P}_\Omega(\mathcal{O}),
\end{aligned}
\label{PTNN_TC}
\end{equation}
where $\mathcal{O},\mathcal{X},\Omega\in\mathbb{R}^{n_1\times n_2\times n_3}$ are respectively the observed data and the underlying recover result, a binary support indicator tensor. Zeros in $\Omega$ indicate the missing entries in the observed tensor. $\mathcal{P}_\Omega(\mathcal{O}) = \Omega\odot\mathcal{O}$  is the elementwise multiplication (Hardamard product) between the
support tensor $\Omega$ and the observed tensor $\mathcal{Y}$. The constraint implies
that the estimated tensor $\mathcal{X}$ agrees with the observed tensor $\mathcal{O}$ in the observed entries.


Let
\begin{equation}
\mathcal{I}_\Phi(\mathbf{\mathcal{X}})=\left\{
\begin{aligned}
&0,\quad       &\text{if}\ \mathbf{\mathcal{X}}\in\Phi,\\
&\infty, &\text{otherwise},
\end{aligned}
\right.
\end{equation}
where $\Phi :=\{\mathbf{\mathcal{X}}\in\mathbb{R}^{n_1\times n_2\times n_3}:\ \mathcal{P}_\Omega(\mathcal{X})=\mathcal{P}_\Omega(\mathcal{O})\}$.
Thus, the tensor completion model in (\ref{PTNN_TC}) can be rewritten as:
\begin{equation}
\min\limits_{\mathbf{\mathcal{X}}} \quad\mathcal{I}_\Phi(\mathbf{\mathcal{X}})+\|\mathbf{\mathcal{X}}\|_{\text{PSTNN}}.\\
\label{PTNN_TC_un}
\end{equation}

After introducing a auxiliary tensor, the problem (\ref{PTNN_TC_un}) is equivalent to
\begin{equation}
\begin{aligned}
\min\limits_{\mathbf{\mathcal{X}}} \quad&\mathcal{I}_\Phi(\mathbf{\mathcal{Y}})+\|\mathbf{\mathcal{X}}\|_{\text{PSTNN}}\\
\text{s.t.}\quad& \mathbf{\mathcal{Y}}=\mathbf{\mathcal{X}}.
\end{aligned}
\label{PTNN_TC_AU}
\end{equation}

The augmented Lagrangian function of (\ref{PTNN_TC_AU}) is given as:
\begin{equation}
\begin{aligned}
L_\beta(\mathbf{\mathcal{X}},\mathbf{\mathcal{Y}},\mathbf{\mathcal{M}})=&
\mathcal{I}_\Phi(\mathbf{\mathcal{Y}})+
 \|\mathbf{\mathcal{X}}\|_{\text{PSTNN}}+
 \langle\mathbf{\mathcal{M}},\mathbf{\mathcal{X}}-\mathbf{\mathcal{Y}}\rangle+
 \frac{\beta}{2}\|\mathbf{\mathcal{X}}-\mathbf{\mathcal{Y}}\|_F^2\\
 =&
 \mathcal{I}_\Phi(\mathbf{\mathcal{Y}})+
 \|\mathbf{\mathcal{X}}\|_{\text{PSTNN}}+
 \frac{\beta}{2}\|\mathbf{\mathcal{X}}-\mathbf{\mathcal{Y}}+\frac{\mathbf{\mathcal{M}}}{\beta}\|_F^2+\mathcal{C},\\
 \end{aligned}
 \label{TC_AUG}
\end{equation}
where $\mathbf{\mathcal{M}}$ is the Lagrangian multiplier, $\beta$ is the Lagrange penalty parameter, and $\mathcal{C} = -\frac{\beta}{2}\|\frac{\mathcal{M}}{\beta}\|_F^2$  is constant with respect to $\mathcal{X}$ and $\mathcal{Y}$. Following the framework of ADMM \cite{boyd2011distributed}, which has shown its effectiveness for solving large scale optimization problems \cite{li2018low,yang2020remote,xie2018kronecker}, we then iteratively update the variables $\mathcal{X}$,  $\mathcal{Y}$ by solving corresponding subproblems and the multiplier $\mathcal{M}$.

\paragraph{Step 1: updating $\mathcal{X}$} The $\mathcal{X}$-subproblem is
\begin{equation}
{\mathbf{\mathcal{X}}}^{k+1} =\arg\left\{\underset{\mathbf{\mathcal{X}}}{\min}\left( \|\mathbf{\mathcal{X}}\|_{\text{PSTNN}}+
 \frac{\beta}{2}\|\mathbf{\mathcal{X}}-{\mathbf{\mathcal{Y}}}^k+\frac{{\mathbf{\mathcal{M}}}^k}{\beta}\|_F^2\right)\right\},
\end{equation}
the solution of which can be exactly calculated by Algorithm \ref{KEY_alg}.

\paragraph{Step 2: updating $\mathcal{Y}$} The $\mathcal{Y}$-subproblem is
\begin{equation}
\mathcal{Y}^{k+1} = \arg\left\{\underset{\mathbf{\mathcal{Y}}}{\min}\left(\mathcal{I}_\Phi(\mathbf{\mathcal{Y}})+
\frac{\beta}{2}\|\mathbf{\mathcal{X}}^{k+1}-\mathbf{\mathcal{Y}}+\frac{\mathcal{M}^k}{\beta}\|_F^2\right)\right\}.
\end{equation}
By minimizing the $\mathcal{Y}$-subproblem, we have $1_\Phi(\mathcal{Y}) = 0$, i.e., $\mathcal{Y} \in \Phi$.  Thus, the solution of the $\mathcal{Y}$-subproblem is given as follows:
\begin{equation}\label{equ:xsolution}
  \begin{cases}
    \begin{aligned}
          &\mathcal{P}_\Omega(\mathcal{Y}^{k+1}) = \mathcal{P}_\Omega(\mathcal{O}),\\
          &\mathcal{P}_{\Omega^{C}}(\mathcal{Y}^{k+1})= \mathcal{P}_{\Omega^{C}}( {\mathbf{\mathcal{X}}}^{k+1}+\frac{\mathcal{M}^{k}}{\beta}),
    \end{aligned}
  \end{cases}
\end{equation}
where $\Omega^{C}$ denotes the complementary set of $\Omega$.

\paragraph{Step 3: updating multiplier} According to the standard ADMM, the multiplier is updated as follows:
\begin{equation}
\mathcal{M}^{k+1}={\mathbf{\mathcal{M}}}^{k}+\beta({\mathbf{\mathcal{X}}}^{k+1}-{\mathbf{\mathcal{Y}}}^{k+1}).
\end{equation}

Finally, Algorithm 2 presents the pseudocode for solving the proposed PSTNN-based tensor completion (TC) model.
\begin{algorithm}[hbtp]
\caption{The pseudocode for solving the PSTNN-based TC model (\ref{PTNN_TC}) by ADMM }
\begin{algorithmic}[1]
\renewcommand{\algorithmicrequire}{\textbf{Input:}} 
\renewcommand{\algorithmicensure}{\textbf{Output:}}
\Require
The observed tensor $\mathbf{\mathcal{O}}\in\mathbb{R}^{n_1\times n_2\times n_3}$, the support of the observed entries $\Omega$, the given tubal multi-rank $\text{rank}_r$, stopping criterion $\epsilon$, the Lagrange penalty parameter $\beta$.
\renewcommand{\algorithmicrequire}{\textbf{Initialization:}}
\Require $\mathbf{\mathcal{X}}^0=\tt{rand}(n_1\times n_2\times n_3)$, $\mathcal{P}_\Omega(\mathcal{X}^0)=\mathcal{P}_\Omega({\mathcal{O}})$, $\mathbf{\mathcal{Y}}^0=\mathbf{\mathcal{X}}^0$, $\mathbf{\mathcal{M}}^0=\tt{zeros}(n_1\times n_2\times n_3)$.

\While {not converged}

\State update $\mathbf{\mathcal{X}}^{k+1}$ with $\left({\mathbf{\mathcal{Y}}}^k-\frac{{\mathbf{\mathcal{M}}}^k}{\beta}\right)$ and $\tau=\frac{n_3}{\beta}$ by algorithm \ref{KEY_alg}

\State ${\mathbf{\mathcal{Y}}}^{k+1} \leftarrow\mathcal{P}_\Omega(\mathcal{O})+\mathcal{P}_{\Omega^{C}}( {\mathbf{\mathcal{X}}}^{k+1}+\frac{\mathcal{M}^{k}}{\beta})$

\State ${\mathbf{\mathcal{M}}}^{k+1} \leftarrow{\mathbf{\mathcal{M}}}^{k}+\beta({\mathbf{\mathcal{X}}}^{k+1}-{\mathbf{\mathcal{Y}}}^{k+1})$


\State Check the convergence conditions $\| \mathbf{\mathcal{X}}^{k+1}-\mathbf{\mathcal{X}}^{k}\|_\infty\leq\epsilon$, $\| \mathbf{\mathcal{Y}}^{k+1} -\mathbf{\mathcal{Y}}^{k} \|_\infty\leq\epsilon$, $\| \mathbf{\mathcal{X}}^{k+1}-\mathbf{\mathcal{Y}}^{k+1}\|_\infty\leq\epsilon$

\EndWhile
\Ensure
The completed tensor $\mathbf{\mathcal{X}}\in\mathbb{R}^{n_1\times n_2\times n_3}$.

\end{algorithmic}
\label{TC_alg}
\end{algorithm}


\subsection{Tensor RPCA using PSTNN}
As mentioned previously, the goal of the tensor RPCA problems is to recover the low-rank tensors from sparsely corrupted observations. A tensor RPCA model using PSTNN can be formulated as
\begin{equation}
\begin{aligned}
\min\limits_{\mathbf{\mathcal{L}},\mathbf{\mathcal{E}}} \quad&\|\mathbf{\mathcal{L}}\|_{\text{PSTNN}}+\lambda\|\mathbf{\mathcal{E}}\|_1 \\
\text{s.t.}\quad& \mathbf{\mathcal{O}}=\mathbf{\mathcal{L}}+\mathbf{\mathcal{E}},
\end{aligned}
\label{PTNN_RPCA}
\end{equation}
where $\mathcal{O},\mathcal{L},\mathcal{E}\in\mathbb{R}^{n_1\times n_2\times n_3}$ are the observed data, the low-rank part, and the sparse corruptions, respectively, and $\lambda$ is a non-negative parameter. Here, we minimize $\|\mathbf{\mathcal{E}}\|_1$, which is the $\ell_1$ norm of $\mathcal{E}$, i.e., the sum of absolute values of entries in $\mathcal{E}$, to enhance the sparsity of $\mathcal{E}$.

The augmented Lagrangian function of \eqref{PTNN_RPCA}is
\begin{equation}
\begin{aligned}
L_\beta(\mathbf{\mathcal{L}},\mathbf{\mathcal{E}},\mathbf{\mathcal{M}})=&
\|\mathbf{\mathcal{L}}\|_{\text{PSTNN}}+\lambda\|\mathbf{\mathcal{E}}\|_1+\langle\mathbf{\mathcal{M}},\mathbf{\mathcal{O}}-\mathbf{\mathcal{L}}-\mathbf{\mathcal{E}}\rangle+
 \frac{\beta}{2}\|\mathbf{\mathcal{O}}-\mathbf{\mathcal{L}}-\mathbf{\mathcal{E}}\|_F^2\\
 =&
\|\mathbf{\mathcal{L}}\|_{\text{PSTNN}}+\lambda\|\mathbf{\mathcal{E}}\|_1+
 \frac{\beta}{2}\|\mathbf{\mathcal{O}}-\mathbf{\mathcal{L}}-\mathbf{\mathcal{E}}-\frac{\mathbf{\mathcal{M}}}{\beta}\|_F^2+\mathcal{C},\\
 \end{aligned}
 \label{RPCA_AUG}
\end{equation}
where $\beta$ is the Lagrange parameter, $\mathbf{\mathcal{M}}$ is the Lagrangian multiplier, and $\mathcal{C} = -\frac{\beta}{2}\|\frac{\mathcal{M}}{\beta}\|_F^2$  is constant with respect to $\mathcal{L}$ and $\mathcal{E}$.


Similar to the updating scheme in the previous section, we then iteratively update the variables $\mathcal{X}$,  $\mathcal{Y}$ by solving corresponding subproblems and the multiplier $\mathcal{M}$, using ADMM \cite{boyd2011distributed}.

\paragraph{Step 1: updating $\mathcal{L}$} The $\mathcal{L}$-subproblem is exhibited as follows:
\begin{equation}
{\mathbf{\mathcal{L}}}^{k+1} =\arg\left\{\underset{\mathbf{\mathcal{L}}}{\min}\left( \|\mathbf{\mathcal{L}}\|_{\text{PSTNN}}+
 \frac{\beta}{2}\|\mathcal{O}-\mathcal{L}-{\mathbf{\mathcal{E}}}^k+\frac{{\mathbf{\mathcal{M}}}^k}{\beta}\|_F^2\right)\right\}.
\end{equation}
Again, we utilize Algorithm \ref{KEY_alg} to solve this subproblem.

\paragraph{Step 2: updating $\mathcal{E}$} The $\mathcal{E}$-subproblem is
\begin{equation}\label{e-subproblem}
\mathcal{E}^{k+1} = \arg\left\{\underset{\mathbf{\mathcal{E}}}{\min}\left(\lambda\|\mathbf{\mathcal{E}}\|_1+
 \frac{\beta}{2}\|\mathbf{\mathcal{O}}-\mathcal{L}^{k+1}-\mathbf{\mathcal{E}}- \frac{\mathcal{M}^k}{\beta}\|_F^2\right)\right\}.
\end{equation}
The solution of \eqref{e-subproblem} can be obtained with the soft-thresholding operator as:
\begin{equation}
\mathcal{E}^{k+1} = \mathbf{\mathcal{S}}_\frac{\lambda}{\beta}\left[ \mathbf{\mathcal{O}}-{\mathbf{\mathcal{L}}}^{k+1}+\frac{{\mathbf{\mathcal{M}}}^{k}}{\beta} \right].
\end{equation}

\paragraph{Step 3: updating multiplier} The multiplier is updated as follows:
\begin{equation}
\mathcal{M}^{k+1}={\mathbf{\mathcal{M}}}^{k}+\beta({\mathbf{\mathcal{O}}}-{\mathbf{\mathcal{L}}}^{k+1}-{\mathbf{\mathcal{E}}}^{k+1}).
\end{equation}

Algorithm \ref{TRPCA_alg} shows the pseudocode for solving the proposed PSTNN-based tensor robust component analysis (TRPCA) model.
\begin{algorithm}[hbtp]
\caption{The pseudocode for solvingthe PSTNN-based TRPCA model (\ref{PTNN_RPCA}) by ADMM }
\begin{algorithmic}[1]
\renewcommand{\algorithmicrequire}{\textbf{Input:}} 
\renewcommand{\algorithmicensure}{\textbf{Output:}}
\Require
The observed tensor $\mathbf{\mathcal{O}}\in\mathbb{R}^{n_1\times n_2\times n_3}$, the given tubal multi-rank $\text{rank}_r$, parameter $\lambda$, stopping criterion $\epsilon$, the Lagrange penalty parameter $\beta$.

\renewcommand{\algorithmicrequire}{\textbf{Initialization:}}
\Require $\mathbf{\mathcal{L}}^0=\mathbf{\mathcal{O}}$, $\mathbf{\mathcal{E}}^0=\mathbf{\mathcal{M}}^0=\tt{zeros}(n_1\times n_2\times n_3)$.

\While {not converged}

\State update $\mathbf{\mathcal{X}}^{k+1}$ with $\left(\mathbf{\mathcal{O}}-\mathbf{\mathcal{L}}-{\mathbf{\mathcal{E}}}^k-\frac{{\mathbf{\mathcal{M}}}^k}{\beta}\right)$ and $\tau=\frac{n_3}{\beta}$ by algorithm \ref{KEY_alg}

\State ${\mathbf{\mathcal{E}}}^{k+1} \leftarrow\mathbf{\mathcal{S}}_\frac{\lambda}{\beta}\left[ \mathbf{\mathcal{O}}-{\mathbf{\mathcal{L}}}^{k+1}+\frac{{\mathbf{\mathcal{M}}}^{k}}{\beta} \right]$

\State ${\mathbf{\mathcal{M}}}^{k+1}\leftarrow{\mathbf{\mathcal{M}}}^{k}+\beta({\mathbf{\mathcal{O}}}-{\mathbf{\mathcal{L}}}^{k+1}-{\mathbf{\mathcal{E}}}^{k+1})$


\State Check the convergence conditions $\| \mathbf{\mathcal{L}}^{k+1}-\mathbf{\mathcal{L}}^{k}\|_\infty\leq\epsilon$, $\| \mathbf{\mathcal{E}}^{k+1} -\mathbf{\mathcal{E}}^{k} \|_\infty\leq\epsilon$, $\| \mathbf{\mathcal{L}}^{k+1}+\mathbf{\mathcal{E}}^{k+1}-\mathbf{\mathcal{O}}\|_\infty\leq\epsilon$

\EndWhile
\Ensure
The low PSTNN tensor $\mathbf{\mathcal{L}}$ and the sparse tensor $\mathbf{\mathcal{E}}$
%
\end{algorithmic}
\label{TRPCA_alg}
\end{algorithm}




\section{Experiments}\label{sec:Exp}

To examine the performance of the proposed methods, we compare the proposed methods \footnote{Our Matlab code is available at \url{https://github.com/TaiXiangJiang/PSTNN}.
} with the TNN-based methods \footnote{Corresponding codes can be downloaded at \url{https://sites.google.com/site/canyilu/} and \url{https://sites.google.com/site/jamiezeminzhang/}.}on the simulated data and different real-world data \footnote{In this paper, we only conduct experiments on the third-order tensors. However, as the t-SVD framework, which is originally suggested for third-order tensors, has been extended for tensors with arbitrary dimensions \cite{martin2013order,zheng2018tensor}. Therefore, the proposed methods can be generalized for high order tensors.}
We adopt two quantitative assessments to accurately measure the quality of the reconstructions.
The first one is the peak signal-to-noise ratio (PSNR), which can be computed by
PSNR is defined as
\[
\text{PSNR}=10\log_{10}\frac{\bar{\mathbf{\mathcal{Y}}}_{\text{GT}}^{2}}{\frac{1}{n^2}\|\mathbf{\mathcal{Y}}-\mathbf{\mathcal{Y}}_{\text{GT}}\|_{F}^{2}},
\]
where $\mathbf{\mathcal{Y}}_{\text{GT}}$, $\bar{\mathbf{\mathcal{Y}}}_{\text{GT}}$, and $\mathbf{\mathcal{Y}}$ are respectively the ground truth tensor, the maximum pixel value of the ground truth tensor, and the reconstructed tensor.
The second one is the structural similarity index (SSIM) \cite{ssim}.
The larger values of PSNR and SSIM are corresponding to the higher quality of the results.
All the numerical experiments are conducted on a PC with a 3.30 GHz CPU and 16 GB RAM.

Throughout our experiments, we assume that the $N$ in \eqref{PTNN_def} is known.
We directly use the ground truth $N$ in the synthetic experiments, while we estimate $N$ by counting the number of the largest
1\% singular values of the first slice of the clean tensor after fast Fourier transform along the third direction.
When the clean data is unavailable, we recommend the heuristic strategies proposed in \cite{wen2012solving}, i.e., the rank-decreasing scheme and the rank-increasing scheme.
These two strategies start to form an overestimated or an underestimated $N$, and then dynamically adjust the estimation by QR decomposition.
The effectiveness of these two strategies has also been validated in \cite{li2019low,Xu2013Tmac}.

\subsection{Synthetic data}

To synthesize the ground-truth tensor, we perform a t-prod $\mathbf{\mathcal A} = \mathbf{\mathcal P}*\mathbf{\mathcal Q}$, where $\mathbf{\mathcal P}\in\mathbb{R}^{n_1\times r\times n_3}$ and $\mathbf{\mathcal Q}\in\mathbb{R}^{r\times n_2\times n_3}$ are independently sampled from an i.i.d. Gaussian distribution $\mathcal{N}(0,\frac{1}{\sqrt{n_1\times n_3}})$. Then, the tubal multi-rank of tensor $\mathbf{\mathcal{A}}\in\mathbb{R}^{n_1\times n_2\times n_3}$ is $[r,r,\cdots,r]^\top$.

\subsubsection{Tensor completion}\label{Sec:sTC}
\begin{figure}[!htp]
  \centering
  \begin{tabular}{cc}
  TNN&PSTNN\\
  \includegraphics[width=0.48\textwidth]{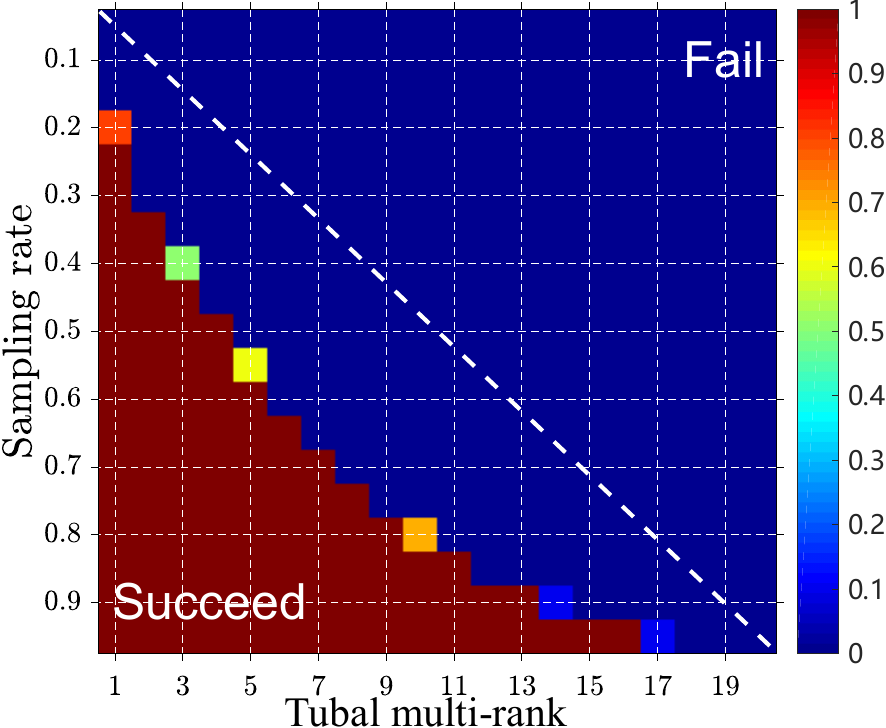}&
  \includegraphics[width=0.48\textwidth]{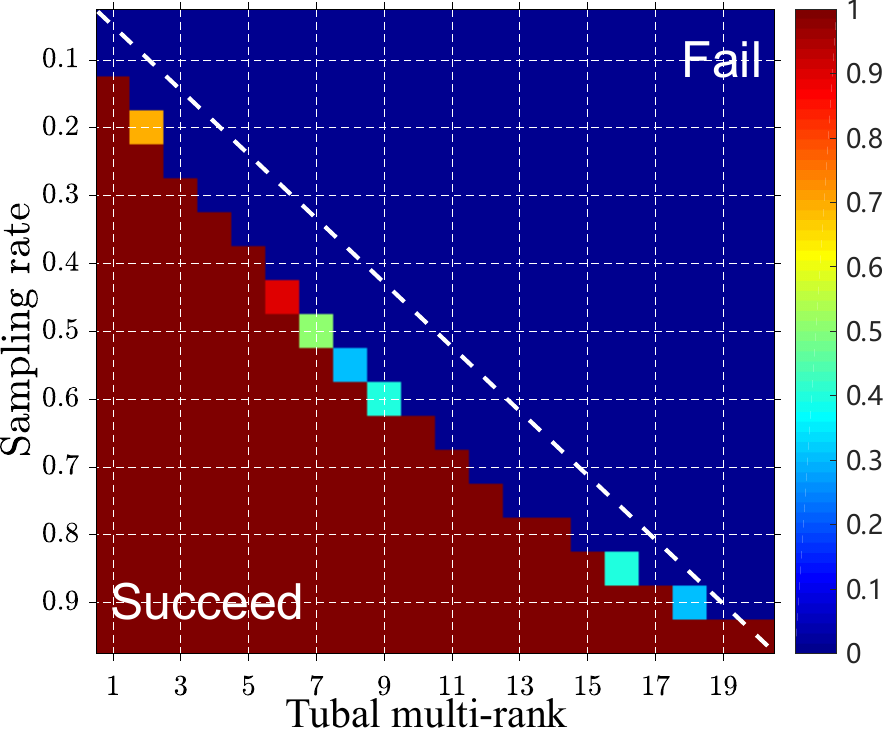}\\
  \multicolumn{2}{c}{(a) $30\times30\times20$ tensor} \\
  TNN&PSTNN\\
  \includegraphics[width=0.48\textwidth]{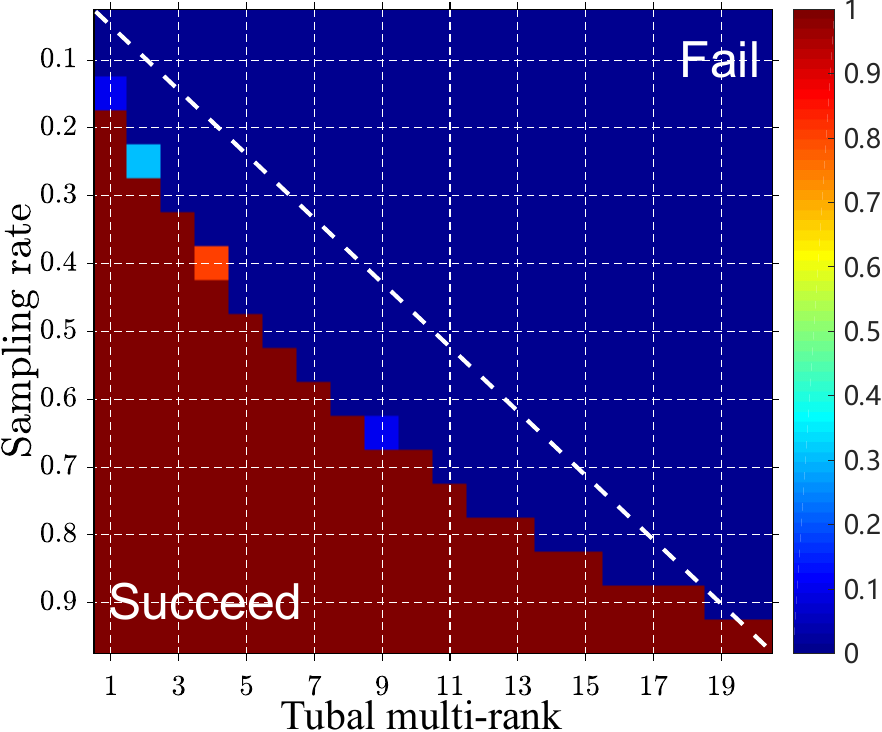}&
  \includegraphics[width=0.48\textwidth]{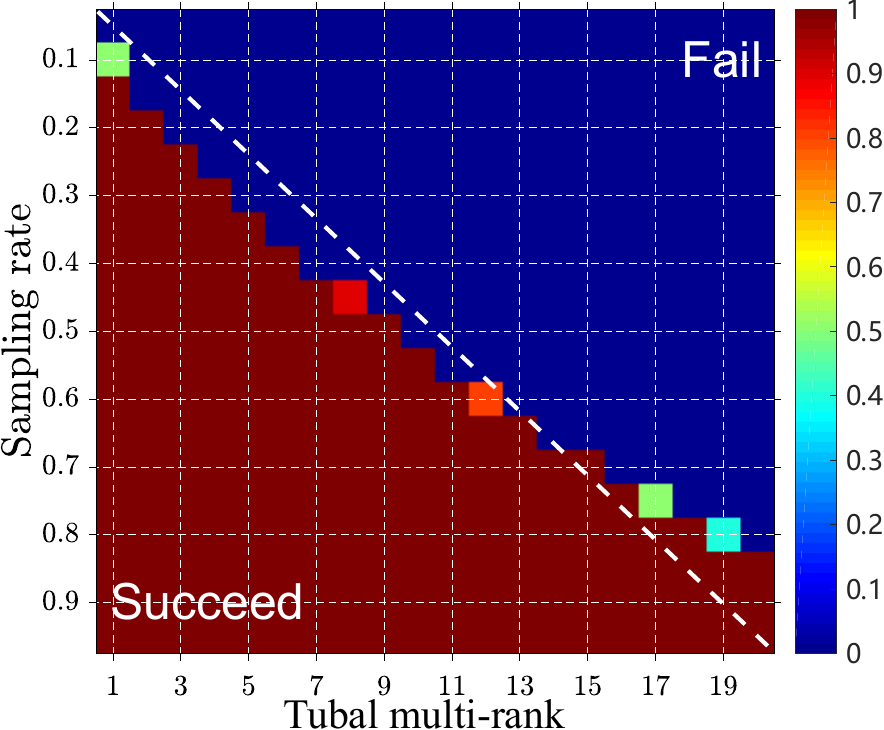}\\
   \multicolumn{2}{c}{ (b) $40\times40\times20$ tensor}
  \end{tabular}
  \caption{Success ratio for synthetic data of two different size and varying tubal multi-ranks with varying sampling rate.
  The left figures illustrate the empirical recovery rate by minimizing the TNN while the right figures by minimizing the PSTNN.
  Each entry in the figures reflects the proportion of the successful recoveries when conducting 10 independent experiments. The white dashed lines are placed on the diagonal line for easier comparison.}
  \label{FTTC}
\end{figure}
For the tensor completion task, we try to recover $\mathcal A$ from the partial observation which is randomly sampled $m$ entries of $\mathcal A$. 
To verify the robustness of the TNN-based TC method and the proposed PSTNN-based TC method, we conducted the experiments with respect to data sizes, the tubal multi-rank $\text{rank}_r$, the sampling rate, {\em i.e.} $\frac{m}{n_1\times n_2\times n_3}$, respectively.
We examine the performance by counting the number of successes.
If the relative square error of the recovered $\widehat{\mathcal A}$ and the ground truth $\mathcal A$, {\em i.e.} $\frac{\|\mathcal{A}-\widehat{\mathcal A}\|^2_F}{\|\mathcal{A}\|^2_F}$, is less than $10^{-3}$, then the recovery is counted as a successful one.
We repeat each case 10 times, and each cell in Figure \ref{FTTC} reflects the success percentage, which is computed by the successful times dividing 10.
Figure \ref{FTTC} illustrates that the proposed PSTNN-based TC method is more robust than the TNN- based TC method, because of bigger brown areas.

\subsubsection{Tensor robust principal components analysis}
\begin{figure}[htp]
  \centering
  \begin{tabular}{cc}
  TNN&PSTNN\\
  \includegraphics[width=0.48\textwidth]{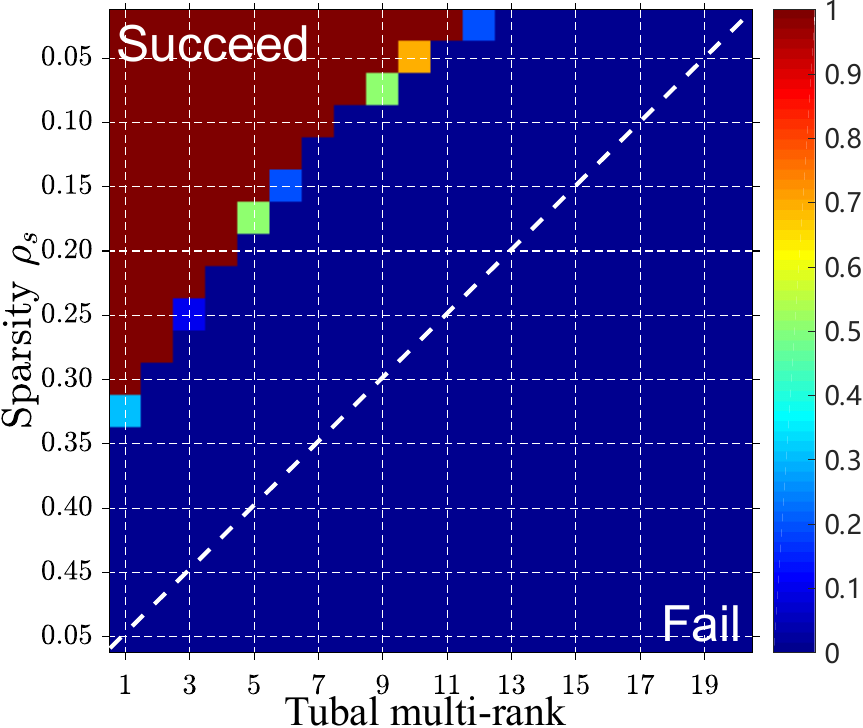}&
  \includegraphics[width=0.48\textwidth]{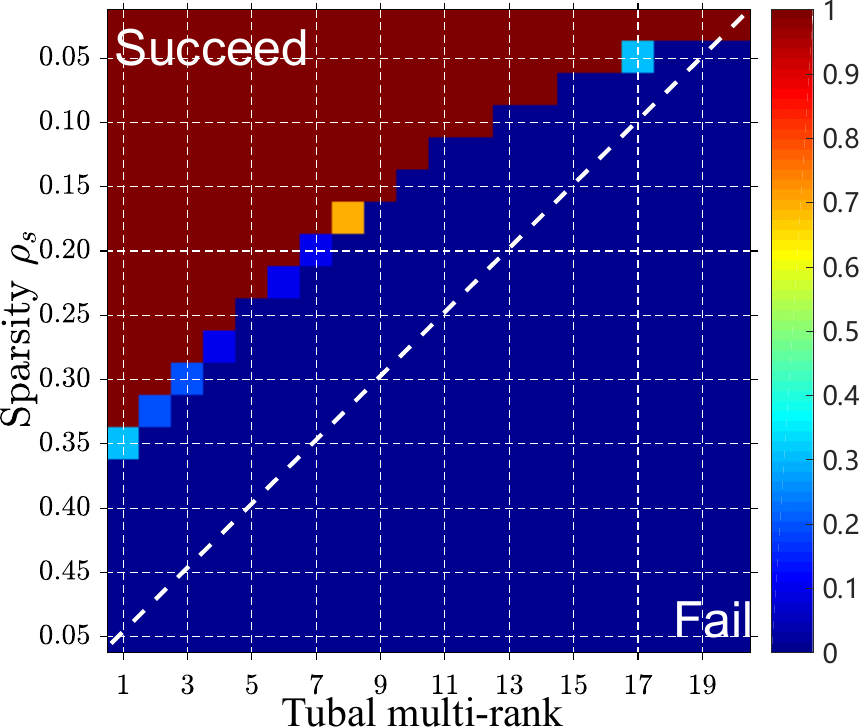}\\
  \multicolumn{2}{c}{(a) $40\times40\times20$ tensor}\\

  TNN&PSTNN\\
  \includegraphics[width=0.48\textwidth]{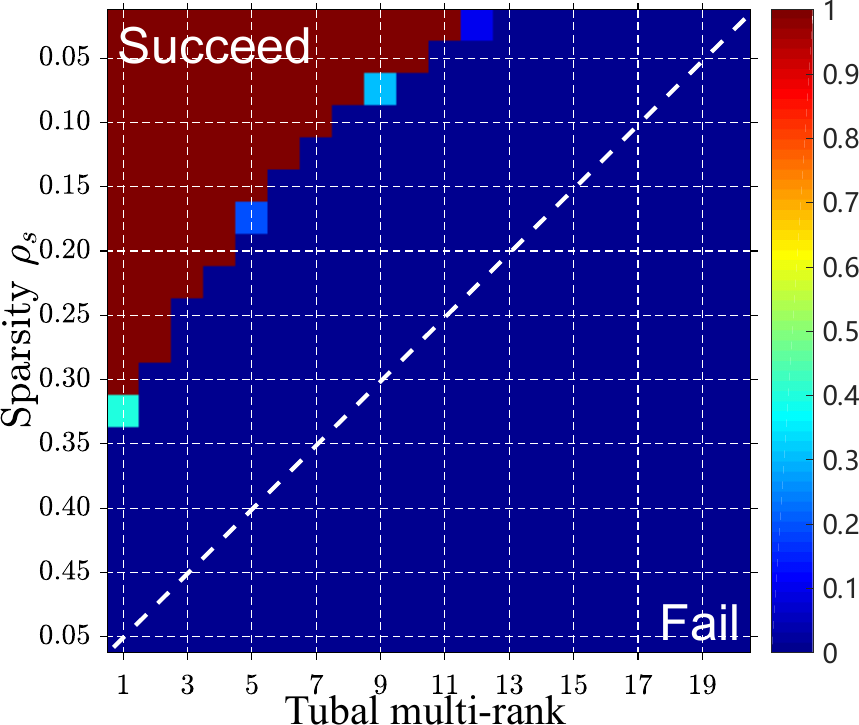}&
  \includegraphics[width=0.48\textwidth]{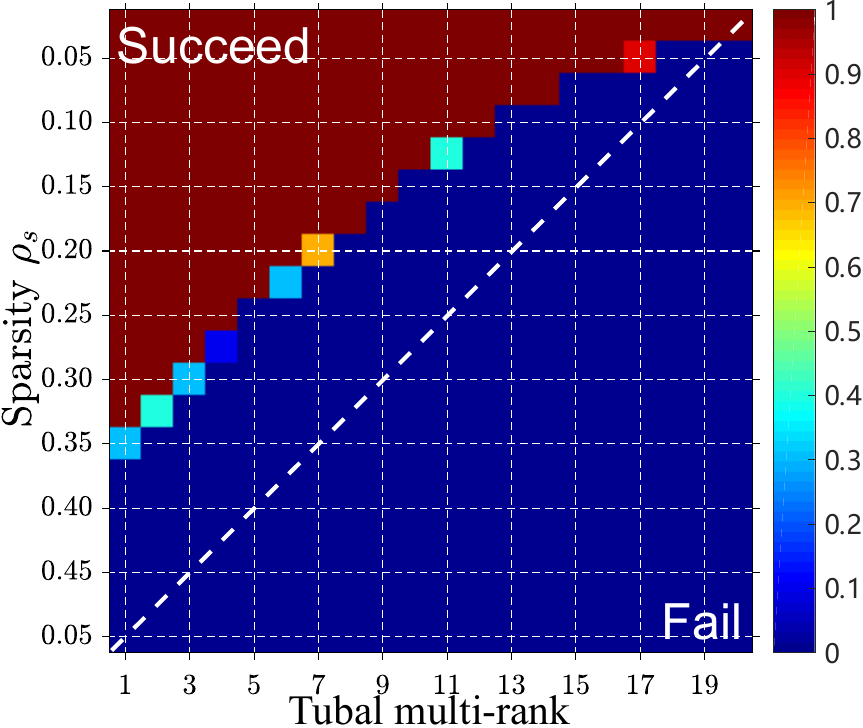}\\
   \multicolumn{2}{c}{ (b) $40\times40\times30$ tensor}
  \end{tabular}
  \caption{The exact recovery results when sparsity and the tubal multi-rank are varying.
  Each entry in the figures reflects the proportion of the successful recoveries when conducting 10 independent experiments.
   The white dashed lines are placed on the diagonal line for easier comparison.}
  \label{FTTRPCA}
\end{figure}

For the tensor robust principal components analysis task, $\mathcal A$ is corrupted by a sparse noise with sparsity $\rho_s$ and uniform distributed values.
We try to recover $\mathcal A$ using Algorithm \ref{TRPCA_alg} and the TNN-based tensor completion method.
The setting of the experiments in this part is similar to that in Section \ref{Sec:sTC}.
We conducted the experiments with respect to data sizes, the tubal multi-rank, sparsity $\rho_s$, respectively.
We report the exact recovery results in Figure \ref{FTTRPCA}.
We repeat each case 10 times, and each cell in Figure \ref{FTTRPCA} reflects the success percentage, which is computed by the successful times dividing 10.
From Figure \ref{FTTRPCA}, we can find that our PSTNN TC method is more robust than the TNN-based TC method, because of the smaller blue areas.

\subsubsection{Sensitivity to initialization}

\begin{figure}[htp]
  \centering
  \includegraphics[width=0.95\textwidth]{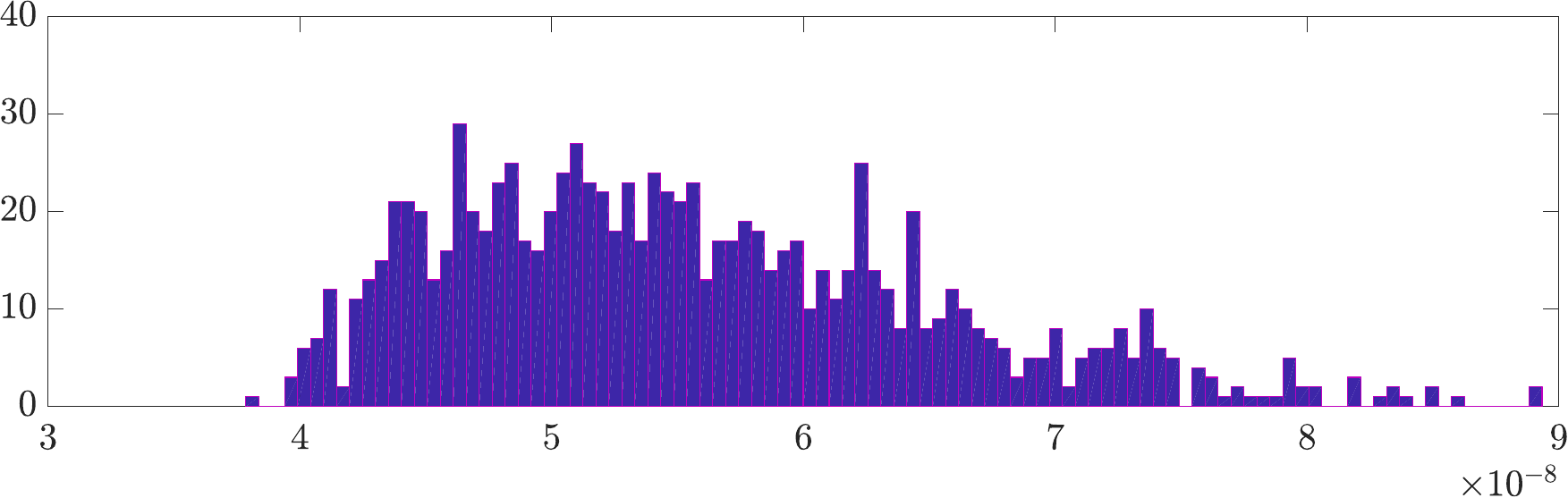}
  \caption{The histogram of  rooted relative squared errors when the initializations are randomly set for the TC task.}
  \label{Initial}
\end{figure}

The converged solution may be different with different initializations, on account of that the proposed objective function is non-convex.
It is necessary to examine the sensitivity of the proposed method against different initializations.
In this subsection, to recover a $25\times 25\times 30$ tensor with tubal multi-rank $\mathbf 5$ and with $10\%$ missing entries in the TC task,
we randomly initialize the tensor for 1000 times. 
The distribution of the rooted relative squared errors is shown in Figure \ref{Initial}.
Although the convergence the proposed algorithms has not been proved with the theoretical guarantee, we can observe from Figure \ref{Initial} that the distribution of the solutions with different initializations concentrates on the near region of the ground truth.

\subsection{Tensor completion for the real-world data}

In this subsection, we conduct experiments on the real-world data, including three video data (``{\it pedestrian}''\footnote{\url{http://www.changedetection.net}}, ``{\it news}'', and ``{\it hall}''\footnote{\url{http://trace.kom.aau.dk/yuv/index.html}}), the MRI data\footnote{\url{http://brainweb.bic.mni.mcgill.ca/brainweb/selection_normal.html}} and the multispectral image (MSI) data\footnote{\url{http://www1.cs.columbia.edu/CAVE/databases/multispectral}}. The compared methods consist of HaLRTC \cite{Liu2013PAMItensor}, the TNN-based TC method \cite{zhang2017exact}, and our PSTNN-based TC method.
The ratio of the missing entries is set as 80\%.
Figure \ref{videoTC} exhibits one frame/band/slice of the completion results.
From Figure \ref{videoTC}, we can conclude that the visual quality of the results obtained by our PSTNN-based TC method is higher than those by HaLRTC and TNN.
The quantitative comparisons are shown in Table \ref{video_comp}, our method obtained the best results with respect to PSNR and SSIM.
The outstanding performance of the PSTNN method on varied real-world data illustrates that the PSTNN is a more precise characterization of the low tubal multi-rank structure.
\begin{table}[htbp]
\setlength{\tabcolsep}{5pt}\renewcommand\arraystretch{0.8}
\centering
\caption{Quantitative comparisons of the completion results by HaLRTC, TNN and PSTNN on the real-world data.}
\begin{tabular}{cccccccc}
\toprule
\multicolumn{2}{c}{Data}  & Size                                        & Index & Observed & HaLRTC  & TNN     & PSTNN       \\ \midrule
\multirow{6}{*}{Video}& \multirow{2}{*}{``{\it pedestrian}''} &\multirow{2}{*}{$158\times 238\times 24$}
                                                     & PSNR  &  7.1475  & 22.6886 & 26.2793 & \bf26.7292 \\
                                                    &&& SSIM  & 0.0459   & 0.6786  & 0.8187  & \bf 0.8288\\ \cmidrule{2-8}
                      & \multirow{2}{*}{``{\it news}''} &\multirow{2}{*}{$158\times 238\times 24$}
                                                     & PSNR  &  9.7447  & 29.9569 & 31.9288 & \bf32.7566 \\
                                                    &&& SSIM  & 0.0618   & 0.9100  & 0.9236  & \bf 0.9296\\ \cmidrule{2-8}
                      & \multirow{2}{*}{``{\it hall}''} &\multirow{2}{*}{$158\times 238\times 24$}
                                                     & PSNR  & 5.5882  & 31.6016 & 33.4691 & \bf34.0856 \\
                                                    &&& SSIM  & 0.0244   & 0.9585  & 0.9605 & \bf 0.9744\\ \midrule

\multicolumn{2}{c}{\multirow{2}{*}{MRI}} &\multirow{2}{*}{$181\times 217\times40$}
                                                     & PSNR  & 10.3162  & 24.3162 & 26.9626 & \bf27.9680 \\
                                                    &&& SSIM  &  0.0887  & 0.7175  & 0.8144  & \bf0.8236\\ \midrule

\multicolumn{2}{c}{\multirow{2}{*}{MSI}} &\multirow{2}{*}{$256\times 256\times31$}
                                                     & PSNR  & 13.8113  & 24.0003 & 28.9523 & \bf30.8586 \\
                                                    &&& SSIM  &  0.1353  & 0.6703  & 0.8695  & \bf0.9026\\ \bottomrule
\end{tabular}
\label{video_comp}
\end{table}

\begin{figure}[htbp]
\scriptsize\setlength{\tabcolsep}{2pt}\renewcommand\arraystretch{0.8}\centering
  \begin{tabular}{ccccc}
Original& Observed & HaLRTC & TNN & PSTNN\\
\includegraphics[width=0.19\textwidth]{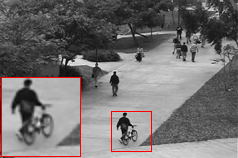}&
\includegraphics[width=0.19\textwidth]{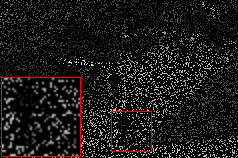}&
\includegraphics[width=0.19\textwidth]{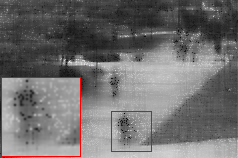}&
\includegraphics[width=0.19\textwidth]{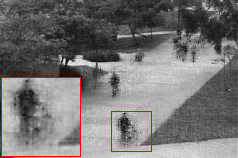}&
\includegraphics[width=0.19\textwidth]{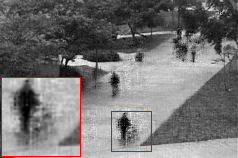}\\

\includegraphics[width=0.19\textwidth]{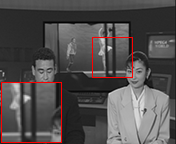}&
\includegraphics[width=0.19\textwidth]{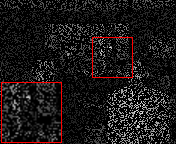}&
\includegraphics[width=0.19\textwidth]{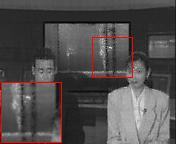}&
\includegraphics[width=0.19\textwidth]{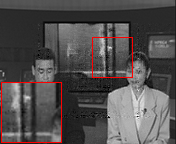}&
\includegraphics[width=0.19\textwidth]{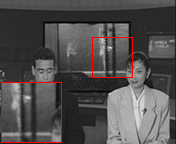}\\

\includegraphics[width=0.19\textwidth]{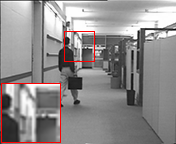}&
\includegraphics[width=0.19\textwidth]{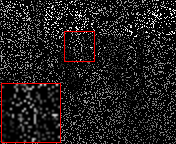}&
\includegraphics[width=0.19\textwidth]{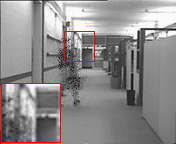}&
\includegraphics[width=0.19\textwidth]{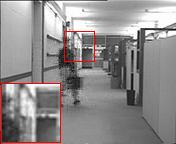}&
\includegraphics[width=0.19\textwidth]{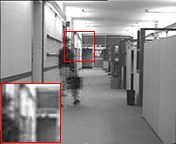}\\

\includegraphics[width=0.19\textwidth]{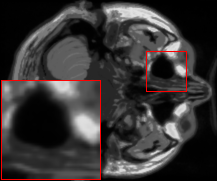}&
\includegraphics[width=0.19\textwidth]{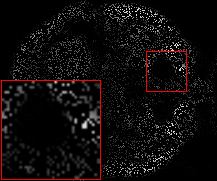}&
\includegraphics[width=0.19\textwidth]{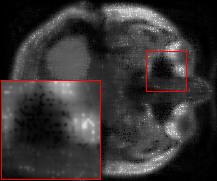}&
\includegraphics[width=0.19\textwidth]{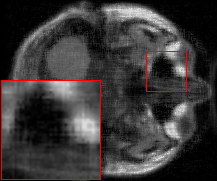}&
\includegraphics[width=0.19\textwidth]{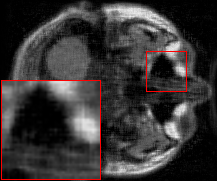}\\

\includegraphics[width=0.19\textwidth]{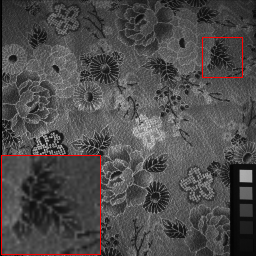}&
\includegraphics[width=0.19\textwidth]{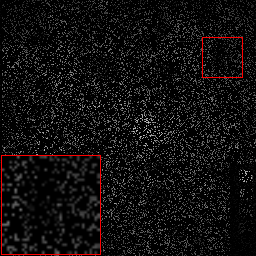}&
\includegraphics[width=0.19\textwidth]{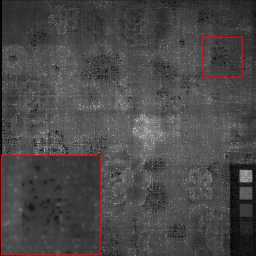}&
\includegraphics[width=0.19\textwidth]{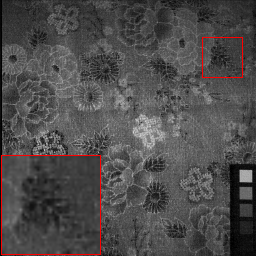}&
\includegraphics[width=0.19\textwidth]{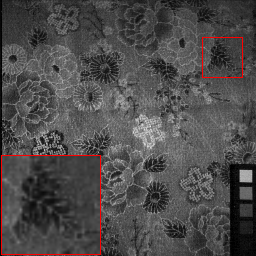}\\

  \end{tabular}
  \caption{Results for the tensor completion for the real-world data. From top to bottom: one frame of the video data (``{\it pedestrian}'', ``{\it news}'', and ``{\it hall}''), one slice of the MRI data, one band of the MSI data.}
  \label{videoTC}
\end{figure}

\subsection{Tensor robust principal components analysis for the color image recovery}

In this subsection, we test the TRPCA methods on the task of the color image recovery. Each image is corrupted by the sparse noise with sparsity 0.2.
We compare our PSTNN-based TRPCA method with the sum of nuclear norm \cite{Liu2013PAMItensor}(SNN)-based TRPCA method and TNN-based TRPCA method \cite{lu2016tensor} on the 4 high quality color images from the Kodak PhotoCD Dataset \footnote{\url{http://r0k.us/graphics/kodak/}} and the homepage\footnote{\url{https://github.com/canyilu/LibADMM}} of the author of \cite{lu2016tensor}.

\begin{table}[htp]
\setlength{\tabcolsep}{8pt}\renewcommand\arraystretch{0.8}
\centering
\caption{Quantitative comparisons of the image recovery results of SNN, TNN and PSTNN on the image data.}
\begin{tabular}{ccccccc}
\toprule
Image &Size                                         & Index & Observed & SNN & TNN & PSTNN\\ \midrule
\multirow{2}{*}{``\it starfish''} &\multirow{2}{*}{$481\times 321\times3$}   & PSNR  & 14.8356 & 25.8286 & 26.411 & \bf28.8492 \\
                                                     && SSIM  & 0.5085 & 0.9440 & 0.9495 & \bf0.9596\\ \midrule
\multirow{2}{*}{``\it door''} &\multirow{2}{*}{$256\times 256\times3$}       & PSNR  & 14.9029 & 27.9449 & 31.4588 & \bf33.4505 \\
                                                     && SSIM  & 0.6200 & 0.9777 & 0.9882 & \bf0.9918\\ \midrule
\multirow{2}{*}{``\it hat1''} &\multirow{2}{*}{$256\times 256\times3$}       & PSNR  & 15.7203 & 23.7104 & 26.3266 & \bf28.5517\\
                                                     && SSIM  & 0.4649 & 0.8993 & 0.9498 & \bf0.9559\\ \midrule
\multirow{2}{*}{``\it hat2''} &\multirow{2}{*}{$256\times 256\times3$}       & PSNR  & 15.3731 & 28.1310 & 31.4964 & \bf32.1626\\
                                                     && SSIM  & 0.4086 & 0.9654 & 0.9789 & \bf0.9810 \\ \bottomrule
\end{tabular}
\label{img_denoising}
\end{table}

\begin{figure}[htbp]
\setlength{\tabcolsep}{2pt}\renewcommand\arraystretch{0.8}
\centering
  \begin{tabular}{ccccc}
Original& Observed & SNN & TNN & PSTNN\\
\includegraphics[width=0.18\textwidth]{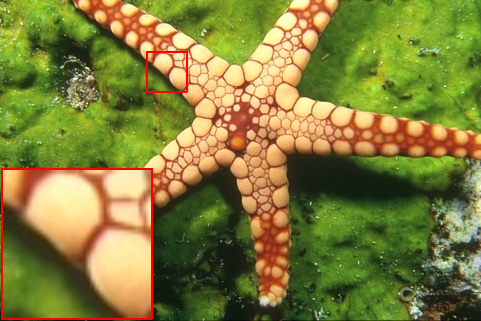}&
\includegraphics[width=0.18\textwidth]{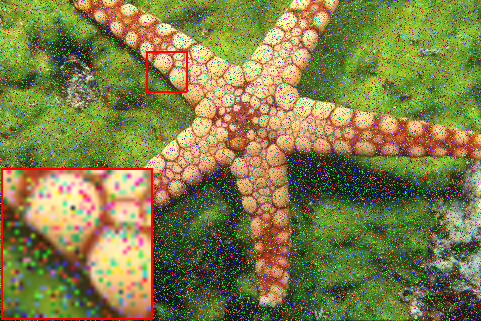}&
\includegraphics[width=0.18\textwidth]{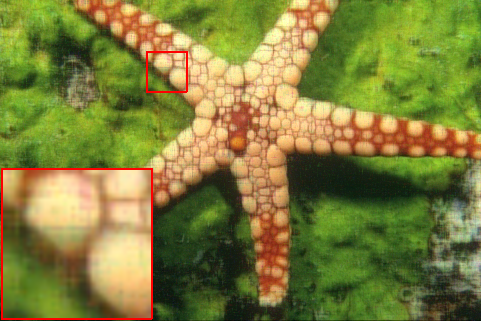}&
\includegraphics[width=0.18\textwidth]{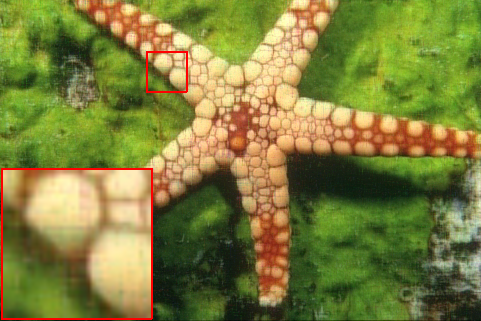}&
\includegraphics[width=0.18\textwidth]{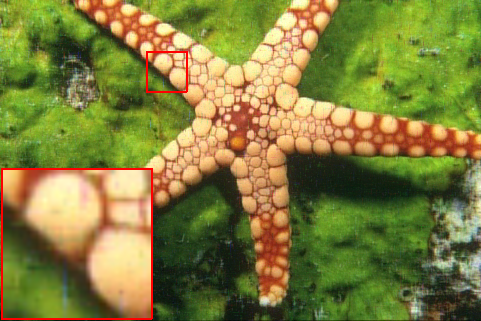}\\
\includegraphics[width=0.18\textwidth]{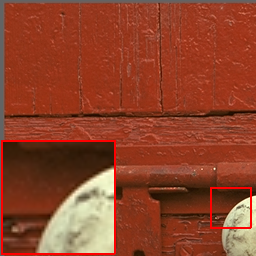}&
\includegraphics[width=0.18\textwidth]{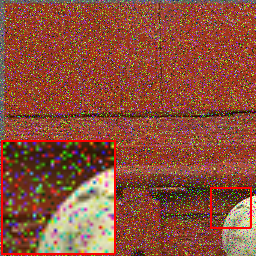}&
\includegraphics[width=0.18\textwidth]{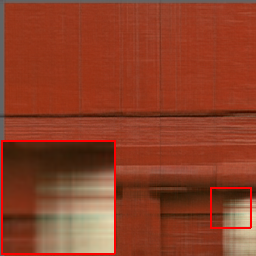}&
\includegraphics[width=0.18\textwidth]{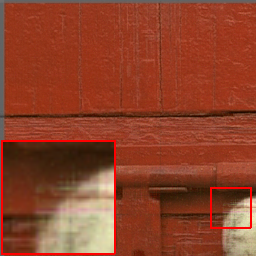}&
\includegraphics[width=0.18\textwidth]{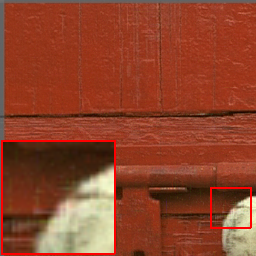}\\
\includegraphics[width=0.18\textwidth]{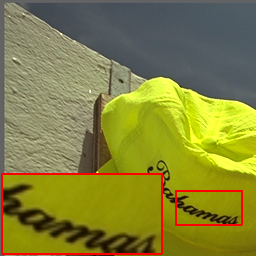}&
\includegraphics[width=0.18\textwidth]{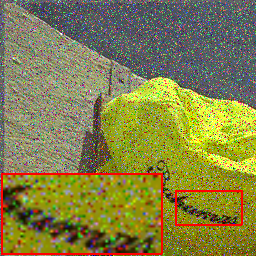}&
\includegraphics[width=0.18\textwidth]{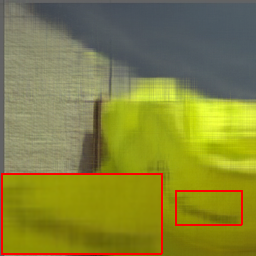}&
\includegraphics[width=0.18\textwidth]{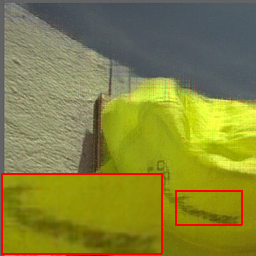}&
\includegraphics[width=0.18\textwidth]{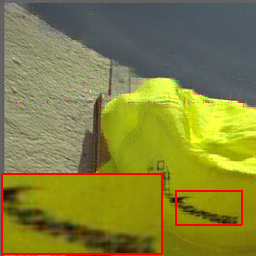}\\
\includegraphics[width=0.18\textwidth]{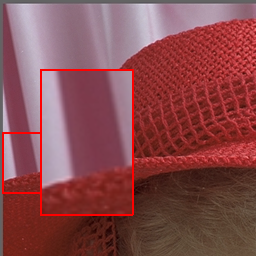}&
\includegraphics[width=0.18\textwidth]{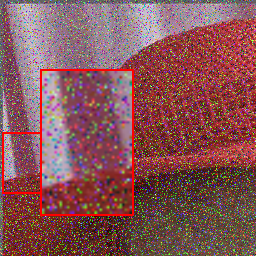}&
\includegraphics[width=0.18\textwidth]{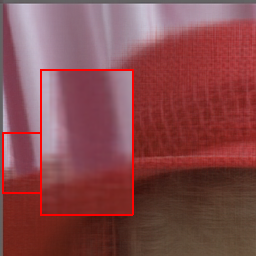}&
\includegraphics[width=0.18\textwidth]{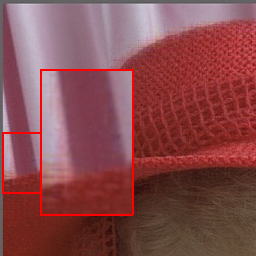}&
\includegraphics[width=0.18\textwidth]{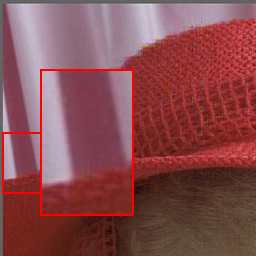}\\

  \end{tabular}
  \caption{Results for the image recovery task.}
  \label{imageTRPCA}
\end{figure}

\begin{figure}[htbp]
  \centering
  \includegraphics[width=0.95\textwidth]{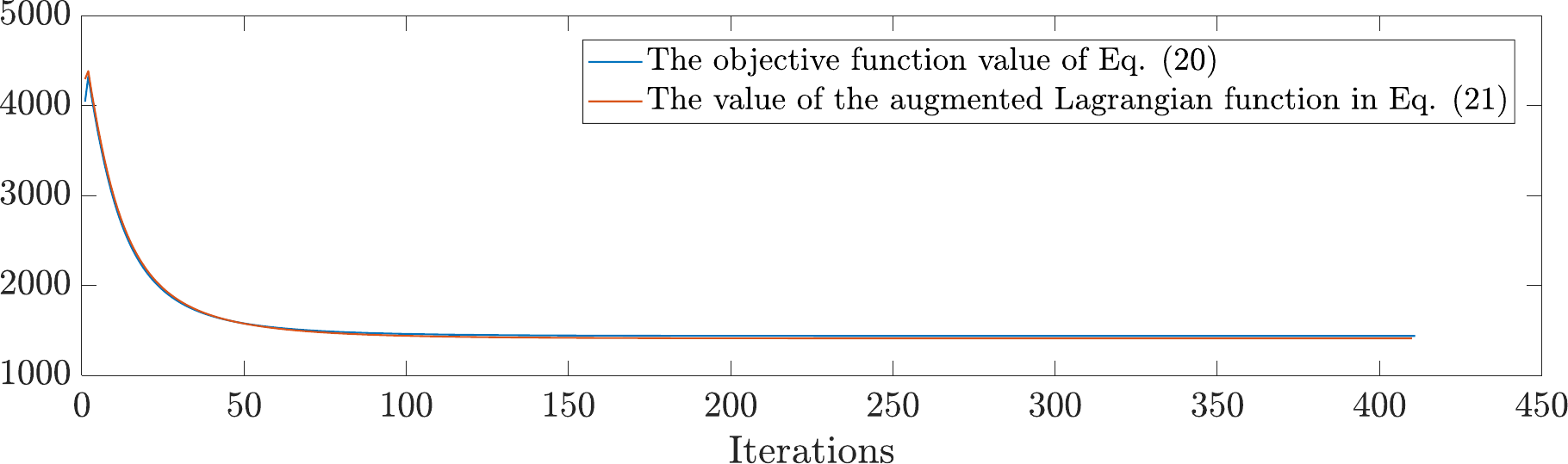}
  \caption{The objective function value in Eq. (20) and the value of the augmented  Lagrangian function in Eq. (21) with respect to iterations when dealing with the image ``{\it starfish}''. (Parameters: $\beta = 1$, $\lambda = 1/\sqrt{\max(n1,n2)*n3}=0.0263$, and $\epsilon = 10^{-7}$.)}
  \label{obj:decay}
\end{figure}

The results are shown in Figure \ref{imageTRPCA}.
From Figure \ref{imageTRPCA}, we can find the results obtained by our PSTNN-based TRPCA method is of higher visual quality, considering the preservation of the image details and textures.
The SNN-based TRPCA method tends to output blurry results.
As for quantitative comparisons exhibited in Table \ref{img_denoising}, our method obtained the best results with respect to PSNR and SSIM while The TNN-based TRPCA method achieves the second-best place. The comparison in this subsection illustrates that our PSTNN-based TRPCA method is more efficient and robust than the SNN-based and TNN-based TRPCA methods.

Meanwhile, as the objective function in Eq. \eqref{PTNN_RPCA} is non-convex, our algorithm based on ADMM should be considered as a local optimization method \cite{boyd2011distributed}. To show its effectiveness of minimizing the objective function, we exhibit the changing of the objective function value, as well as the value of the augmented Lagrangian function in Eq. \eqref{RPCA_AUG}, with respect to iterations in Figure \ref{obj:decay}. The decaying curves illustrate that our algorithm effectively minimizes the objective function value.

\section{Conclusions}\label{sec:Con}
In this paper we propose a novel surrogate of the tensor tubal multi-rank, {\em i.e.}, PSTNN, within the t-SVD framework.
We extend the PSVT operator for the matrices in the complex field to solve the proposed PSTNN-based minimization problem, which is fundamental for solving the subsequent PSTNN-based tensor recovery models.
Two PSTNN-based minimization models for tensor completion and tensor robust principal component analysis are proposed.
Two efficient ADMM algorithms, using the PSVT solver, have been developed to solve the models.
The effectiveness of the proposed PSTNN-based methods is illustrated by the experiments on the data of various types.

\section*{Acknowledgments}
The authors would like to thank the editor and reviewers for giving us many comments and suggestions, which are of great value for improving the quality of this manuscript.
This work is supported in part by the National Natural Science Foundation of China under Grant 61772003, Grant 61876203 and Grant 61702083, in part by the Fundamental Research Funds for the Central Universities under Grant  ZYGX2016J132, Grant ZYGX2016J129 and Grant ZYGX2\-016KYQD142, and in part by the Science Strength Promotion Programme of UESTC.

\section{Appendix}
\subsection{The definitions in the t-SVD framework} \label{tsvddefs}
\begin{mydef}[t-product \cite{kilmer2011factorization}]
The t-product $\mathbf{\mathcal{C}}=\mathbf{\mathcal{A}}*\mathbf{\mathcal{B}}$
of $\mathbf{\mathcal{A}}\in \mathbb{R}^{n_{1}\times n_2\times n_{3}}$ and $\mathbf{\mathcal{B}}\in \mathbb{R}^{n_{2}\times n_4\times n_{3}}$ is a tensor of size
$n_1\times n_4 \times n_3$, where the $(i,j)$-th tube $\mathbf{c}_{ij:}$ is given by
\begin{equation}
\mathbf{c}_{ij:} = \mathbf{\mathcal{C}}(i,j,:) = \sum_{k=1}^{n_2}\mathbf{\mathcal{A}}(i,k,:)\circledast\mathbf{\mathcal{B}}(k,j,:)
\end{equation}
where $\circledast$ denotes the circular convolution between two tubes of same size.
\end{mydef}


Interpreted in another way, a 3-D tensor of size $n_1\times n_2 \times n_3$
can be viewed as a $n_1\times n_2$ matrix with treating the basic units as a tube.
In the t-prod of two tensors, the interaction of the two basic units is the circular convolution instead of the multiplication.

\begin{mydef}[tensor conjugate transpose \cite{kilmer2011factorization}]
The conjugate transpose of a tensor $\mathbf{\mathcal{A}}\in \mathbb{R}^{n_{2}\times n_1\times n_{3}}$ is tensor $\mathbf{\mathcal{A}}^\text{\rm H}\in \mathbb{R}^{n_{1}\times n_2\times n_{3}}$ obtained by conjugate transposing each of the frontal slice and then reversing the order of transposed frontal slices 2 through $n_3$:
\begin{equation*}
\begin{aligned}
\left(\mathbf{\mathcal{A}}^\text{\rm H}\right)^{(1)}&=\left(\mathbf{\mathcal{A}}^{(1)}\right)^\text{\rm H}\quad\text{and}\\
\left(\mathbf{\mathcal{A}}^\text{\rm H}\right)^{(i)}&=\left(\mathbf{\mathcal{A}}^{(n_3+2-i)}\right)^\text{\rm H},\quad i=2,\cdots,n_3.
\end{aligned}
\end{equation*}
\end{mydef}


\begin{mydef}[identity tensor \cite{kilmer2011factorization}]
The identity tensor $\mathbf{\mathcal{I}}\in \mathbb{R}^{n_{1}\times n_1\times n_{3}}$ is defined as a tensor whose first frontal slice is the $n_1\times n_1$ identity matrix, and the other frontal slices are zero matrices.
\end{mydef}


\begin{mydef}[orthogonal tensor \cite{kilmer2011factorization}]
A tensor $\mathbf{\mathcal{Q}} \in \mathbb{R}^ {n_{1} \times n_1\times n_{3}}$  is an orthogonal tensor if
\begin{equation}
\mathbf{\mathcal{Q}}^\text{\rm H}*\mathbf{\mathcal{Q}}=\mathbf{\mathcal{Q}}*\mathbf{\mathcal{Q}}^\text{\rm H}=\mathbf{\mathcal{I}}.
\end{equation}
\end{mydef}


\begin{mydef}[block diagonal form \cite{kilmer2011factorization}]\label{Def:bldg}
$\overline{\mathcal{{A}}}$ is used to denote the block-diagonal form unfolding of the Fourier transformed tensor of $\mathcal{A}$, i.e., $\widehat{\mathbf{\mathcal{A}}}$. That is
\begin{equation}
\begin{aligned}
\overline{\mathcal{{A}}}&\triangleq \tt{blockdiag}(\widehat{\mathbf{\mathcal{A}}})\\
&\triangleq
\left [
\begin{tabular}{cccc}
$\widehat{\mathbf{\mathcal{A}}}^{(1)}$&&&\\
&$\widehat{\mathbf{\mathcal{A}}}^{(2)}$ &&\\
&&$\ddots$ &\\
&&&$\widehat{\mathbf{\mathcal{A}}}^{(n_3)}$
\end{tabular}\right]
\in\mathbb{C}^{n_1n_3\times n_2n_3}.
\end{aligned}
\end{equation}
\end{mydef}

It is not difficult to find that $\overline{\mathcal{{A}}^\text{\rm H}}=\overline{\mathcal{{A}}}^\text{\rm H}$, i.e.,  the block diagonal form of a tensor's conjugate transpose equals to the matrix conjugate transpose of the tensor's block diagonal form. Further more, for any tensor $\mathbf{\mathcal{A}}\in \mathbb{R}^{n_{1}\times n_2\times n_{3}}$ and $\mathbf{\mathcal{B}}\in \mathbb{R}^{n_{2}\times n_4\times n_{3}}$, we have
\[
\mathbf{\mathcal{A}}*\mathbf{\mathcal{B}}=\mathbf{\mathcal{C}} \Leftrightarrow \overline{\mathcal{A}}\cdot\overline{\mathcal{B}}=\overline{\mathcal{{C}}},
\]
where $\cdot$ is the matrix product.


\begin{mydef}[f-diagonal tensor \cite{kilmer2011factorization}]
We call a tensor $\mathbf{\mathcal{A}}\in \mathbb{R}^{n_{1}\times n_2\times n_{3}}$ f-diagonal if all of its frontal slices are the diagonal matrices.
\end{mydef}


\begin{thm}[t-SVD \cite{kilmer2011factorization}]
For $\mathbf{\mathcal{A}}\in \mathbb{R}^{n_{1}\times n_2\times n_{3}}$, the t-SVD of $\mathbf{\mathcal{A}}$ is as the following form
\begin{equation}
\mathbf{\mathcal{A}}=\mathbf{\mathcal{U}}*\mathbf{\mathcal{S}}*\mathbf{\mathcal{V}}^\text{\rm H}
\end{equation}
where $\mathbf{\mathcal{U}}\in \mathbb{R}^{n_{1}\times n_1\times n_{3}}$ and $\mathbf{\mathcal{V}}\in \mathbb{R}^{n_{2}\times n_2\times n_{3}}$ are both orthogonal tensors, and $\mathbf{\mathcal{S}}\in \mathbb{R}^{n_{1}\times n_2\times n_{3}}$ is an f-diagonal tensor.
\end{thm}

\begin{figure}[!htp]
  \centering
  \includegraphics[width=0.65\textwidth]{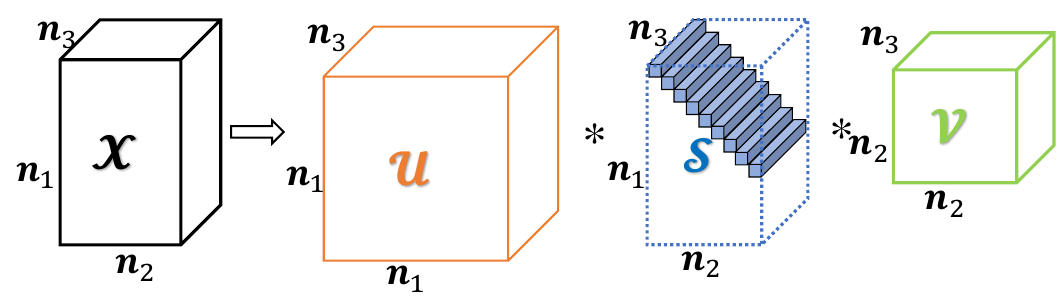}
  \caption{The t-SVD of an $n_1 \times n_2 \times n_3$ tensor.}
  \label{tsvd}
\end{figure}

The t-SVD is illustrated in Figure \ref{tsvd} and can be efficiently computed by the frontal slice wise singular value decomposition (SVD) after Fourier transform.

\begin{mydef}[tubal multi-rank \cite{Oguz2014Tensor}]\label{Def:tubal}
For a three way tensor $\mathbf{\mathcal{A}}\in\mathbb{R}^{n_1\times n_2\times n_3}$, 
its tubal multi-rank, denoted as $\text{rank}_r(\mathbf{\mathcal{A}})$, is defined as a vector, whose $i$-th ($i = 1,2,\cdots,n_3$) element represents the rank of the $i$-th frontal slice of $\widehat{\mathbf{\mathcal{A}}}$ i.e.,
\begin{equation}
\text{rank}_r(\mathbf{\mathcal{A}})=[\text{rank}(\widehat{\mathbf{\mathcal{A}}}^{(1)}), \text{rank}(\widehat{\mathbf{\mathcal{A}}}^{(2)}),\cdots, \text{rank}(\widehat{\mathbf{\mathcal{A}}}^{(n_3)})]^\top.
\end{equation}
\end{mydef}



\begin{mydef}[tubal nuclear norm (TNN) \cite{Oguz2014Tensor}]\label{Def:TNN}
The tubal nuclear norm of a tensor $\mathbf{\mathcal{A}}\in \mathbb{R}^{n_{1}\times n_2\times n_{3}}$, denoted as $\|\mathbf{\mathcal{A}}\|_{\text{\rm TNN}}$, is defined as the sum of singular values of all the frontal slices of $\overline{\mathbf{\mathcal{A}}}$.
\end{mydef}
In particular,
\begin{equation}
\begin{aligned}
\|\mathbf{\mathcal{A}}\|_{\text{TNN}}\triangleq\|\overline{\mathcal{{A}}}\|_{*}=\sum\limits_{i=1}^{n_3}\|\widehat{\mathbf{\mathcal{A}}}^{(i)}\|_*.
\end{aligned}
\label{tnn}
\end{equation}

\subsection{The proof to Theorem \ref{PSVT_TH}}\label{ProofTH3.1}
{\textbf{Proof to Theorem \ref{PSVT_TH}}\it\ \
Lets consider $\mathbf A=\mathbf U_A\mathbf D_A\mathbf V_A^H = \sum_{i=1}^{l} \sigma_i(\mathbf A)\mathbf u_i\mathbf v_i^H$, where $\mathbf U_A=[\mathbf u_1,\cdots,\mathbf u_m]\in \mathcal(U)_m$, $\mathbf V_A=[\mathbf v_1,\cdots,\mathbf v_m]\in \mathcal(V)_n$ and $\mathbf D_A=\text{diag}(\sigma(\mathbf A))$, where the singular values $\sigma(\cdot)=[\sigma_1(\cdot),\cdots,\sigma_l(\cdot)]\geqslant0$ are sorted in a non-increasing order.
Also we define the function $J(\mathbf A)$ as the objective function of \eqref{pssv_matrix}.
The first term of \eqref{pssv_matrix} can be derived as follows:
\begin{equation}\label{First_term}
\begin{aligned}
&\frac{1}{2}\|\mathbf A-\mathbf B\|_F^2=\frac{1}{2}\left(\|\mathbf B\|_F^2-2<\mathbf A,\mathbf B>+\|\mathbf A\|_F^2\right)\\
&=\frac{1}{2}\left(\|\mathbf B\|_F^2-2\sum\limits^l_{i=1}\sigma_i(\mathbf A)\mathbf u_i^H\mathbf B \mathbf v_i +\sum\limits^l_{i=1}\sigma_i(\mathbf A)^2\right)
\end{aligned}
\end{equation}
In the minimization of (\ref{First_term}) with respect to $\mathbf A$, $\|\mathbf B\|_F^2$ is regarded as a constant and thus can be ignored.
For a more detailed representation, we change the parameterization of $\mathbf A$ to $(\mathbf U_A,\mathbf V_A,\mathbf D_A)$ and minimize the function:
\begin{equation}
\begin{aligned}
J(\mathbf U_A,\mathbf V_A,\mathbf D_A) =\frac{1}{2}\sum\limits_{i=1}^l\left(-2\sigma_i(\mathbf A)\mathbf u_i^H\mathbf B \mathbf v_i + \sigma_i(\mathbf A)^2\right)+\tau\sum\limits^l_{i=N+1}\sigma_i(\mathbf A)
\end{aligned}\label{Jx}
\end{equation}

From von Neumann's lemma, the upper bound of $\mathbf u_i^H\mathbf B\mathbf v_i$ is given as $\sigma_i(\mathbf B)=\max\{\mathbf u_i^H\mathbf B\mathbf v_i\}$ for all $i$ when $\mathbf U_A=\mathbf U_B$ and $\mathbf V_A=\mathbf V_B$.
Then (\ref{Jx}) becomes a function depending only on $\mathbf D_A$ as follows:
\begin{equation}
\begin{aligned}
J&(\mathbf U_B,\mathbf V_B,\mathbf D_A) =\frac{1}{2}\sum\limits_{i=1}^l\left(-2\sigma_i(\mathbf A)\sigma_i(\mathbf B) + \sigma_i(\mathbf A)^2\right)+\tau\sum\limits^l_{i=N+1}\sigma_i(\mathbf A)\\
&=\frac{1}{2}\sum\limits_{i=1}^N \left(-2\sigma_i(\mathbf A)\sigma_i(\mathbf B)+ \sigma_i(\mathbf A)^2\right)+\frac{1}{2}\sum\limits_{i=N+1}^l\left(-2\sigma_i(\mathbf A)\sigma_i(\mathbf B) + \sigma_i(\mathbf A)^2+2\tau\sigma_i(\mathbf A)\right).
\end{aligned}\label{Jy}
\end{equation}


Since (\ref{Jy}) consists of simple quadratic equations for each $\sigma_i(\mathbf A)$ independently, it is trivial to show that the minimum of (\ref{Jy}) is obtained at $\hat{\mathbf D}_A=\text{diag}\left(\hat\sigma(\mathbf A)\right)$ by derivative in a feasible domain as the first-order optimality condition, where $\hat\sigma(\mathbf A)$ is defined as
\begin{equation}
\hat\sigma(\mathbf A)=\left\{
\begin{aligned}
&\sigma_i(\mathbf B),\quad       &\quad\text{if} \quad i<N+1,\\
&\max\left(\sigma_i(\mathbf B)-\tau,0\right), &\text{otherwise}.
\end{aligned}
\right.
\end{equation}

Hence, the solution of \eqref{pssv_matrix} is $\mathbf A^*=\mathbf U_B\hat{\mathbf D}_A \mathbf V_B^H$.
This result exactly corresponds to the PSVT operator where a feasible solution $\mathbf A^*=\mathbf U_B(\mathbf D_{B_1}+\mathbf{\mathcal{S}}_\tau[\mathbf D_{B_2}])\mathbf V_B^H$ exists.
$\Box$}\\

\bibliography{ref}
\bibliographystyle{elsarticle-num}

\end{document}